\documentclass{amsart}
\usepackage{amsthm,amsmath,amsfonts,amssymb,amscd,mathrsfs,graphics}
\usepackage{latexsym}
\usepackage{graphicx}
\usepackage{subfigure}

\usepackage{epic}

\usepackage{enumerate}

\newtheorem{thm}{Theorem}[section]

\newtheorem{lem}[thm]{Lemma}
\newtheorem{prop}[thm]{Proposition}
\newtheorem{cor}[thm]{Corollary}

\makeatletter
\renewcommand{\@seccntformat}[1]{\S{\csname
the#1\endcsname}\hspace{0.5em}}
\makeatother

\begin{document}

\title [Schur rings over  ${\bf {\rm Sp}(n,2)}$ ] {Schur rings over  ${\bf {\rm Sp}(n,2)}$ and multiplicity one subgroups}

\author{ Stephen Humphries  }
  \address{Department of Mathematics,  Brigham Young University, Provo, 
UT 84602, U.S.A.
E-mail: steve@mathematics.byu.edu, }
\date{}

\begin{abstract}  
We study commutative Schur rings over the symplectic groups Sp$(n,2)$  containing the class $\mathcal C$  of symplectic transvections.
We find the possible partitions of $\mathcal C$ determined by the Schur ring. We show how this restricts the possibilities for multiplicity one subgroups of Sp$(n,2)$.
\medskip

\noindent {\bf Keywords}: Symplectic group; characteristic $2$; Schur ring; Orthogonal group; symmetric group; $3$-transposition group. \newline 
\medskip
\subjclass[2010]{Primary 05B10. Secondary: 20C05.}
%\subjclass[2010]{Primary 20C15. Secondary: 20E45.}
\end{abstract}

\maketitle

\theoremstyle{plain}

\theoremstyle{definition}
\newtheorem*{dfn}{Definition}
\newtheorem{exa}[thm]{Example}
\newtheorem{rem}[thm]{Remark}

\newcommand{\ds}{\displaystyle}
\newcommand{\bs}{\boldsymbol}
\newcommand{\mb}{\mathbb}
\newcommand{\mc}{\mathcal}
\newcommand{\mf}{\mathfrak}
\renewcommand{\mod}{\operatorname{mod}}
\newcommand{\mult}{\operatorname{Mult}}

\def \a{\alpha} \def \b{\beta} \def \d{\delta} \def \e{\varepsilon} \def \g{\gamma} \def \k{\kappa} \def \l{\lambda} \def \s{\sigma} \def \t{\theta} \def \z{\zeta}

\numberwithin{equation}{section}

\setlength{\leftmargini}{1.em} \setlength{\leftmarginii}{1.em}
\renewcommand{\labelenumi}{\setlength{\labelwidth}{\leftmargin}
   \addtolength{\labelwidth}{-\labelsep}
   \hbox to \labelwidth{\theenumi.\hfill}}

\maketitle

\section{Introduction}

Let $G$ be a finite group. For $X\subseteq G$ we let $\overline {X}=\sum_{x \in X} x \in \mathbb QG$ and $X^{(-1)}=\{x^{-1}:x \in X\}$. Then  a {\it Schur ring} \cite {sch,wie,cur,mp} over $G$ is a subalgebra $\mathfrak S$ of $\mathbb QG$ determined by a partition $C_1=\{1\},C_2,\cdots,C_r$ of $G$ where
(i)  for  $1\le i\le r$ there is $1\le i'\le r$ with $C_i^{(-1)}=C_{i'}$;
(ii) there are   $\lambda_{ijk}\in \mathbb Z^{\ge 0}$ with $\overline {C_i}\,\overline {C_j}=\sum_{k=1}^r \lambda_{ijk}\overline {C_k}.$

The $C_i$ are called {\it principal sets} (or {\it p-sets}) and an {\it $\mathfrak S$-set} is a union of principal sets.  One way to construct a Schur ring over $G$ is to take a subgroup $H\le G$ and let the $C_i$ be the {\it $H$-classes} $g^H=\{g^h:h \in H\}$. This Schur ring is denoted $\mathfrak S(G,H)$. 

With this definition we can propose the following: 

\noindent {\bf Problem}: {\it  Suppose that we have a finite group $G$ that is generated by a conjugacy class $\mathcal C$. Determine all partitions of $\mathcal C$ that come from a commutative Schur ring $\mathfrak S$ over $G$ that contains $ {\mathcal C}$ as an $\mathfrak S$-set. }

This problem is particularly relevant for groups containing a canonical class $\mathcal C$ of generators, and more specifically when those generators are involutions.
Some  examples of such a class of groups are the 3-transposition groups: a group $G$ that has a conjugacy class  $\mathcal C$   of involutary generators satisfying: for $a,b \in \mathcal C$  the order of $ab$ is $1,2$ or $3$ \cite {Fi1971,Asch,Hall}. The most-studied class of 3-transposition groups are the symmetric groups $S_n$ where $\mathcal C$ is the set of transpositions. Other examples of 3-transposition groups are the symplectic groups Sp$(n,2)$ and the orthogonal groups ${\rm SO}^{\pm}(n,2)\le {\rm Sp}(n,2)$.  In \cite {hum} we addressed the above problem for $S_n$ obtaining:

              \begin{thm} \label{thm992} Let $n \ge 6$ and 
   let 
$\mathfrak S$ be a commutative  S-ring over $S_n$ where   $  {\mathcal C}=(1,2)^{S_n}$ is an $\mathfrak S$-set. 
Suppose that $C_1,\dots,C_r$ are the  principal sets of $\mathfrak S$ that are contained in $\mathcal C$. Then $r\le 3$.
  \end{thm}
  
More information  about the partition of $(1,2)^{S_n}$ in the $S_n$ case  is given in \cite {hum}.

In this paper we consider the 3-transposition group $G=G_n:={\rm Sp}(n,2)={\rm Sp}(V)$ where $V$ is a non-degenerate symplectic space over $\mathbb F_2=\{0,1\}$ of even dimension $n\ge 2$. Here the conjugacy class $\mathcal C$ is the class of symplectic transvections $\mathcal T=\mathcal T_n$ that determines ${\rm Sp}(n,2)$ as a 3-transposition group. 

The  main result in this paper   is
              \begin{thm} \label{mainthm1} Let $n \ge 14$ and 
   let 
$\mathfrak S$ be a commutative  S-ring over {\rm Sp}$(n,2)$ where  $ {\mathcal T}$ is an $\mathfrak S$-set. 
Suppose that $C_1,\dots,C_r$ are the  principal sets of $\mathfrak S$ that are contained in $\mathcal T$. Then   
$r=1$ and $C_1=\mathcal T$.
  \end{thm}
  
%  More details regarding  cases (ii), (iii)  of this result are given in Theorems \ref {thmOpos} and \ref {thmOneg}.

We also consider the cases $n=2,4,6$. For a group $G$ and $H \le G$ we say that $(G,H)$ is a {\it strong Gelfand pair} if for every $\chi \in \hat G, \psi \in \hat H$ we have
$\langle \psi,\chi_H\rangle \le 1$. We then say that $H$ is a { \it strong Gelfand subgroup of $G$}. Here $\hat G$ is the set of irreducible characters of $G$. Note that $G$ is always a strong Gelfand subgroup of $G$, and if $H\le G$ is one, then so is $H^g, g \in G$.
We note that the concept of a strong Gelfand pair is   what is sometimes referred to as the `multiplicity one' property \cite{A,AA}.

 Now  by  \cite [Corollary 1']  {Trav} or \cite [Corollary 2.2]  {Kar}   $(G,H)$ is a strong Gelfand pair exactly when $\mathfrak S (G,H)$ is commutative, and in this case, if we have a 3-transposition group $G$ with conjugacy class $\mathcal C$, then $\mathcal C$ will be an $\mathcal S$-set of $\mathfrak S(G,H)$. We   show that for $n=2,4$ the only partitions of $\mathcal T$ are those determined by strong Gelfand pairs.  This result is included in \cite {hum} since  $G_2={\rm Sp}(2,2)\cong S_3, G_4={\rm Sp}(4,2)\cong S_6$:

\begin {thm} \label{thmsp4}  
(A) There are  three  strong Gelfand subgroups of $G_2\cong S_3$, namely $\langle (1,2)\rangle, \langle (1,2,3)\rangle, S_3$.
 Their sizes and the  sizes of their p-sets in $\mathcal T$  are
 
\noindent (i) $2, [1,2]$;
  (ii) $3, [3]$;  (iii) $6, [3]$.

\noindent (B) There are  seven  strong Gelfand subgroups of $G_4\cong S_6,$ each of which is either $G_4$ or a maximal subgroup.
 Their sizes and the  sizes of their p-sets in $\mathcal T$  are

\noindent (i) $48, [3,12]$; (ii) $48, [1,6,8]$;  (iii) $72, [6,9]$; (vi) $120, [5,10]$;
 (v) $120, [15]$;
 (vi) $360, [15]$;
 (vii) $720, [15]$.

\end{thm}

A Magma  \cite {ma}  calculation allows us to find the only strong Gelfand subgroups of $G_6, G_8,G_{10},G_{12}$:

\begin {thm} \label{thmsp633}  There are two  strong Gelfand subgroups of $G_6$.
 Their sizes and the  sizes of their p-sets in $\mathcal T$  are

\noindent (i) $25920, [27,36]$;
\noindent (ii) $1451520, [63]$.

For $n=8,10,12$ there are no proper strong Gelfand subgroups of $G_n$. \qed 

\end{thm}

To prove this it suffices to check the maximal subgroups of each such $G_n$, since if we have $H\le K\le G$ where $H$ is a strong Gelfand subgroup of $G$, then $K$ is also a strong Gelfand subgroup of $G$. 
In general, our main theorem gives:
\begin{cor}\label{cormain} 
Let $H$ be a proper strong Gelfand subgroup of $G_n, n\ge 8$. Then 

\noindent (i) $H$ acts transitively on $\mathcal T$; and

\noindent (ii) $H$ does not contain any symplectic  transvections.
\end{cor}
\noindent{\it Proof} (i)
By Theorem \ref {thmsp633} we can assume $n\ge 14$.
 If $H$ does not act transitively on $\mathcal T$, then the commutative Schur ring $\mathfrak S(G_n,H)$ partitions $\mathcal T$ into at least two subsets, contrary to Theorem \ref {mainthm1}. 
 (ii) If $H$ contains a transvection, then by (i) it will contain all transvections. But then $H=G_n$, a contradiction. \qed\medskip

Here, as indicated above, we could restrict ourselves to the situation where $(G,H)$ is a strong Gelfand pair and $H$ is maximal. Then conditions (i) and (ii) of Corollary \ref {cormain} are satisfied by the maximal subgroups $H<G_n$ studied by Dye in \cite {D1,D2}. These only occur when $n=2m$ with $m$ prime, and they are the stabilizers of a set $\{U_i\}_i$ (called a spread) of subspaces where $\{U_i \setminus \{0\}\}_i$ partitions $V\setminus \{0\}$. We denote this $H$ by ${\rm Dye}_n$ and we show later that  ${\rm Dye}_n$ is not a strong Gelfand subgroup.

In \S\S 2-4 we give definitions relating to the symplectic structure of $V=\mathbb F_2^n$  and results regarding products of symplectic transvections. We also define functions $f_1(C),\cdots,f_4(C)$ depending on a set $C\subset \mathcal T$ and  with domain $C$. We show that these functions are constant on p-sets of our Schur ring. We also show that if $f_2(C)=0$, then $\langle C\rangle $ acts on $C$ by conjugation. 
Then in \S 5 we give an overview of the remainder of the  proof, showing how we reduce to the case $r=1$.

\medskip 
\noindent {\bf Acknowledgement} All calculations made in the preparation of this paper were accomplished using Magma \cite {ma}.

\section {Definitions and Preliminary results}

If $V=\mathbb F_2^n$ is a non-degenerate symplectic space with symplectic form $\cdot$ then there is a (standard) symplectic basis $e_1,\cdots,e_n$, where 
$e_i \cdot e_j=1$ if and only if $i+j=n+1$. 
We will also sometimes write a symplectic basis as $a_1=e_1,b_1=e_n,a_2=e_2,b_2=e_{n-1},\cdots,a_{n/2}=e_{n/2},b_{n/2}=e_{n/2+1}$.
We let  $V^\sharp:=V \setminus \{0\}$.   Then $G_n={\rm Sp}(n,2)$ is the group of linear maps of $V$ that preserve this form \cite {Art,Die}.

To each  $v \in V^\sharp$ we associate the {\it symplectic  transvection} $t_v\in G_n$, where 
$$
(w)t_v=w+(v\cdot w)v \text { for } w \in V. 
$$
If $v=0 \in V$, then we naturally define $t_0=id_G$. Let $\mathcal T=\mathcal T_n$ denote the set of symplectic  transvections in $G_n$. Then $\mathcal T$ is a conjugacy class of involutions in $G$.
For $1\le i,j\le n$ we let $t_i=t_{e_i}, t_{i,j}=t_{e_i+e_j}$ etc. 

In what follows we will    identify $v \in V^\sharp$ with $t_v \in \mathcal T$  and so  identify any $U \subseteq V^\sharp$ with 
$\mathcal T(U):=\{t_u:u \in U\}$. 

%For  $S \subseteq V$ we let $\mathcal T(S)=\{t_s:s \in S \}$, and for $H \le G$  let $\mathcal T(H)=\mathcal T \cap H$.

We recall that a subset $X\subset V$ is {\it isotropic} if $u \cdot v=0$ for all $u,v \in X.$ 
For  $X \subseteq V$ we let $X^\perp=\{u \in V:u\cdot x=0 \text { for all } x \in X \}$.
For $X,Y \subseteq V$ we write $X\perp Y$ if $x \cdot y=0$ for all $x \in X, y \in Y$. 

We will have need to use the {\it Schur-Wielandt principle}: if $\alpha=\sum_{g \in G} f(g) g\in \mathfrak S,$ where $f:G \to \mathbb R$ is the coefficient function, then for each
$\l \in \mathbb R$ we have $\sum_{g \in G:f(g)=\l}g \in \mathfrak S$. 
We will also use the following without reference  \cite {Art,Die}:

\begin{thm}[Witt's Theorem]\label{thmwitt}
Let $U,W$ be subspaces of  $V$ and  let $\varphi:U \to W$ be  an isometry. Then $\varphi$ can be extended to an isometry $\psi:V\to V.$
\end{thm}

We have the following fundamental properties of transvections:: 
\begin{lem}\label{lem111}
Let $t_a,t_b \in \mathcal T$. Then

\noindent (i) $t_a=t_b$ if and only if $a=b$ if and only if the order of $t_at_b$ is $1$.

\noindent (ii) For $a \ne b$ we have:  $a\cdot b=0$ if and only if the order of $t_at_b$ is $2$ if and only if $t_a$ and $t_b$ commute.

\noindent (iii) $a\cdot b=1$ if and only if the order of $t_at_b$ is $3$ if and only if $t_a$ and $t_b$ do not commute if and only if $t_a^{t_b}=t_{a+b}=t_b^{t_a}$.

\noindent (iv) For $\alpha \in G, \a \in V^\sharp,$ we have $\a^{-1}t_a\a=t_{(a)\a}$. 

\noindent (v) For   pairs $t_a,t_b\in \mathcal T, a \ne b,$ and $ t_c,t_d\in \mathcal T,c \ne d,$   the elements $t_at_b$ and $t_ct_d$ are   conjugate if and only if $a\cdot b=c \cdot d$.

\noindent (vi) If $\a=t_{v_1}t_{v_2}\cdots t_{v_s}, v \in V$, then $(v)\a=v+w$, where $w \in {\rm Span} (v_1,\cdots,v_s)$.
\qed
\end{lem}

 \begin{lem}\label{lem1}
 Let $a,b \in V^\sharp$ with $a\cdot b=1.$
 If $t_c,t_d \in \mathcal T$, then $t_ct_d=t_at_b$, if and only if  we have  one of:
 
 \noindent (i) $c=a,d=b$; 
 
\noindent (ii) $c=b,d=a+b$; or 
 
\noindent (iii) $c=a+b,d=a$.
 
 Here   the  pair $[t_c,t_d]$ has the form $[t_a^{g^i},t_b^{g^i}],$ for some $ i=0,1,2,$ where $g=t_at_b$. 
 \end {lem}
 \noindent {\it Proof}
 Note that $t_a^{t_b}=t_{a+b}$. For $c,d$ as in (i)-(iii) above,   one checks that $t_ct_d=t_at_b$.
 
Conversely: assume   $t_ct_d=t_at_b$. If $c \cdot d=0$, then $t_at_b$ and $t_ct_d$ are not conjugate and so not equal.
Thus  $c \cdot d=1$,
and from $(x)t_at_b=(x)t_ct_d, x\in V,$ we get
 \begin{align}
\label {eq1a}
 (b\cdot x)b+(a\cdot x)[a+b] = (d \cdot x) d+ (c \cdot x)[c+d].
 \end{align}

 Putting $x=a$ and $  x=b$ in  Eq. (\ref {eq1a}) we get
 $
 b,a+b \in {\rm {\rm Span}} (c,d)
 $ (respectively),
 so  that ${\rm Span}(a,b)$= {\rm Span}$(c,d)$ and we write:
 $$
 c=\alpha a+\beta b,\qquad d=\gamma a+\delta b \text { where } \a,\b,\g,\d \in \mathbb F_2.
 $$
 
 Putting  $x=a$ in Eq. (\ref {eq1a}) we get
 $$
 b=\delta(\gamma a+\delta b)+\b (\a a+\b b+\gamma a+\d b)
 $$
 Since $a \cdot b=1$,  $a,b$ are linearly independent, and  equating coefficients of $a,b$ we get
 \begin{align}\label{eq231}
 0=\d\gamma +\b\a +\b\g,\qquad 
 1=\d+\g+\d\g.
 \end{align}
 Similarly, putting $x=b$ in    Eq. (\ref {eq1a}),   we get
 \begin{align}\label{eq232}
 1=\g+\a+\a\g,\qquad
 1=\g\d+\a\b+\a\d.
 \end{align}
 Solving Eqs (\ref {eq231}), (\ref  {eq232}) over $\mathbb F_2$ gives
 $[c,d]\in \{[a,b],[b,a+b],[a+b,a]\}$. \qed
 \medskip

 \begin{lem} \label{lem2}
(i) If $a,b \in V^\sharp$ with $a\cdot b=0, a \ne b$, then the only $c,d \in V^\sharp$ with $t_ct_d=t_at_b$ are  where $[c,d]\in \{[a,b], [b,a]\}$.
 
\noindent (ii) For  $a,b,c \in V^\sharp$ we cannot have $t_at_b=t_c$. 

 \end{lem}
 \noindent {\it Proof} (i) This is similar to the   proof of Lemma \ref {lem1}.
 (ii) If $a=b$, then $t_c=1$, a contradiction. If $a \ne b$, then one  shows that the image of $t_at_b-1$ has dimension $2$, a contradiction, since the image of $t_c-1$ has dimension $1$.
 \qed\medskip
 
 The next result tells us  when $\overline C$ and $\overline D$ commute, for disjoint $C, D \subseteq \mathcal T$.

  \begin{lem}\label{lem3}
 For disjoint and non-empty   $C,D \subseteq \mathcal T$ we have: $\overline {C}, \overline {D}$ commute if and only if 
 for all $t_a \in C, t_b \in D$ with $a\cdot b=1$ we have $t_{a+b} \in C \cup D$.
 \end{lem}
  \noindent {\it Proof} First note that if $C \perp D$, then  $\overline {C}, \overline {D}$ commute, so we may now assume that 
   $\overline {C}, \overline {D}$ commute and   there are 
  $t_a \in C, t_b \in D$ with $a\cdot b=1$.  Then $t_at_b \in \overline {C}\cdot  \overline {D}$, and so  there are $t_c \in D, t_d \in C$ with 
  $t_ct_d=t_at_b$. By Lemma \ref {lem1} we have one of:
  
\noindent  (i)  $c=a,d=b$; 
  
\noindent  (ii) $c=b,d=a+b$; or
  
\noindent   (iii) $c=a+b,d=a$.
  
  If  (i), then $t_a \in C \cap D=\emptyset$, a contradiction; if  (ii) or (iii), then  $t_{a+b} \in C \cup D$.
  
   For the converse: assume that if $t_a \in C, t_b \in D, a\cdot b=1$, then  $t_{a+b}=t_a^{t_b}=t_b^{t_a} \in C \cup D$.
We show that $\overline {C}, \overline {D}$ commute. 
 Let $t_a \in C, t_b \in D$ so that $t_at_b \in \overline C\,\overline D$. We first show that $t_at_b$ is a summand of  $\overline {D}\, \overline {C}$. 
 If $a\cdot b=0$, then $t_at_b=t_bt_a\in \overline D\cdot \overline C$ and we are done.
 So assume   $a \cdot b=1$.  
Then $t_a^{t_b}=t_{a+b}  \in C \cup D$. We have two cases:
 
\noindent (a) If $t_a^{t_b} \in C$, then $\overline D \, \overline {C}$ contains $t_bt_a^{t_b}=t_at_b$ as a summand, and we are done.
 
\noindent  (b) If $t_a^{t_b} \in D$, then  $\overline D \, \overline {C}$ contains $t_a^{t_b} t_a=t_b^{t_a}t_a=t_at_b$, as required. 
 
 From Lemma \ref {lem1}, Lemma \ref {lem2} and the above enumeration of cases we see that, for fixed $t_a\in C,t_b\in D,$  there is only one way to find $t_c \in D, t_d \in C$ with $t_at_b=t_ct_d$. Thus the coefficient of $t_at_b$ in  $\overline C\,\overline D$ and in $\overline D\,\overline C$  is $1$, and so $ \overline C\, \overline D= \overline D\, \overline C$.\qed
 \medskip

 \begin{cor} \label{corsubres}
 Let $C,D \subset \mathcal T$ be disjoint and assume that $\overline C,\,\overline D$ commute. Let $W\le V$ be a subspace. Then 
 $\overline {C\cap W},\overline {D\cap W}$ also commute.
 \end{cor} 
 \noindent {\it Proof} We use  Lemma \ref {lem3}: assume that $t_a \in C \cap \mathcal T(W), t_b \in D \cap  \mathcal T(W), a\cdot b=1.$ Then
 $a,b \in W$
  and by Lemma \ref {lem3} applied to the pair
 $t_a\in C, t_b\in D$, we have $t_{a+b} \in C\cup D$. 
 But $W$ is a subspace and so $a+b \in W$; so either
 $t_{a+b}\in C \cap  \mathcal T(W)$ or $t_{a+b}\in D \cap  \mathcal T(W)$.
 \qed 
 
 \medskip

 \begin{lem}\label{lem31}
 Let $C \subseteq \mathcal T$. Choose $a,b \in V^\sharp, a \cdot b =1$. 
 Then $t_at_b, t_bt_a$ occur in $\overline {C}^2$ with the same non-zero coefficient   as follows:
 
 \begin{enumerate}

\item  the coefficient is zero if and only if $|C \cap \{t_a,t_b,t_{a+b}\}|\le 1.$
 
\item   the coefficient is $1$ if and only if $|C \cap \{t_a,t_b,t_{a+b}\}| =2.$

\item   the coefficient is $3$ if and only if $|C \cap \{t_a,t_b,t_{a+b}\}| =3.$
  
  \end{enumerate} 
 \end{lem}
  \noindent {\it Proof}
 By Lemma \ref {lem1} we need only check the cases where $C$ is one of
 $\{t_a\},
 \{t_a,t_b\},$ $ 
 \{t_a,t_{a+b}\},
  \{t_b,t_{a+b}\},
   \{t_a,t_b,t_{a+b}\}$.  Checking each case gives the result.
   \qed
  \medskip 
 
 \begin{cor} \label{corab}
 \noindent (i) Let $a,b \in V^\sharp, a \ne b, a \cdot b=0$. Then the  coefficient of $t_at_b$ in $\overline{\mathcal T}^2$ is $2$.
 
 \noindent (ii) Let $a,b \in V^\sharp, a \ne b, a \cdot b=1$. Then the  coefficient of $t_at_b$ in $\overline{\mathcal T}^2$ is $3$.
 
 In particular, $\overline{\mathcal T}^2=|\mathcal T|+2(t_1t_2)^G+3(t_1t_n)^G$ and    $(t_1t_2)^G, (t_1t_n)^G\in \mathfrak S$.\qed 
   
 \end{cor}

\noindent {\bf Definition}   Let $\overline {C}\in \mathfrak S,C\subseteq \mathcal T$. Define  ${\rm Tr}^2(C)$ to be the set of triples $\{a,b,a+b\}, a,b \in V^\sharp, a \cdot b=1,$  that have
exactly two of $t_a,t_b,t_{a+b}$ in $C$. Let ${\rm Tr}^3(C)$ denote the set 
  of all $\{a,b,a+b\}, a,b \in V^\sharp, a \cdot b=1$ where 
    $t_a,t_b,t_{a+b}\in C$.
   
   Each  set $U=\{a,b,a+b\}, a \cdot b=1,$  (or $U=\{a,b\}, a \cdot b=1,$)  determines a sum:
$$
   U^* := t_at_b + t_bt_a   =t_at_{a+b}+t_{a+b}t_a=t_bt_{a+b}+t_{a+b}t_b \in \mathbb Q G_n.
   $$ Note that $U^*$ is well-defined.
   The point of this definition is that if $a\cdot b=1$, then
   $$
   (t_a+t_b)^2=2+\{t_a,t_b\}^*, \text { and }    (t_a+t_b+t_{a+b})^2=3+3\{t_a,t_b,t_{a+b}\}^*.
   $$
  Further, if $C\subseteq \mathcal T$, then the only way that $\{t_a,t_b\}^*$ can occur in $\overline {C}^2$ is if 
  $C$ contains at least two of $\{t_a,t_b,t_{a+b}\}$.
We have:

   \begin{lem}\label{lem32}
    Let $\overline {C}\in \mathfrak S,C\subseteq \mathcal T$. Then
    
    \noindent (i) $\sum_{T \in {\rm Tr}^2(C)} T^* \in \mathfrak S$;
    
    \noindent (ii)   $\sum_{T \in {\rm Tr}^3(C)} T^* \in \mathfrak S$.
    
    \end{lem} 
  
   \noindent {\it Proof}  By Corollary \ref {corab} $(t_1t_n)^G \in \mathfrak S$. 
By Lemma \ref {lem31} we have
  $$
  \overline {C}^2=\left ( \sum_{T \in {\rm Tr}^2(C)} T^*\right )+3\left ( \sum_{T \in {\rm Tr}^3(C)} T^*\right )+ Z,
  $$
  where no summand of $Z$ is in $(t_1t_n)^G$. Now  Corollary \ref {corab}  and  the Schur-Wielandt principle give the result.
  \qed
   \medskip

 \begin{lem} \label{lem43} Let $t_a,t_x,t_y,t_{z} \in  \mathcal T$, where $t_xt_yt_{z} = t_a$. Then 
$[x,y,z]$ is one of:
 
\noindent (i) $[a,d,d]$ for some $d \in V^\sharp$;
 
\noindent (ii)  $[d,a,d]$ for some $d \in V^\sharp$ with $a\cdot d=0$;
 
\noindent (iii)  $[d,d,a]$ for some $d \in V^\sharp$;
 
\noindent (iv) $[b,a,a+b]$  for some $b \in V^\sharp$  with $a \cdot b=1$;
 
%\noindent (v) $[a+b,a,b]$  where $a \cdot b=1$;
 
\noindent (v) $[b,a+b,b]$ for some $b \in V^\sharp$ with  $a \cdot b=1$;
 
% \noindent (vii) $[a+b,b,a+b]$  where $a \cdot b=1$;
 
 \end{lem}
\noindent{\it  Proof} 
First, it is easy to check that each of (i)-(v) gives a solution to $t_xt_yt_z=t_a$.

If $x=y$, then $z=a$ and we have (iii). So assume $x \ne y$.
If $y=z$, then $x=a$ and we have (i), so assume $y \ne z$.
If $x=z$, then $t_a=t_xt_yt_x$. 
If $x \cdot y=0$, then  $t_a=t_xt_yt_x=t_y$ and so $y=a$, so we have (ii);
 if $x\cdot y=1$, then
$t_a=t_xt_yt_x=t_{x+y}$, so $a=x+y$, and we have case (v). Now we  assume that $x,y,z$ are distinct.

If $x \cdot y=0$, then $t_xt_y=t_at_z$ and  by Lemma \ref {lem2} $\{x,y\}=\{a,z\}$. %Thus $a \in \{x,y\}$.
If $x=a$, then $y=z$, a contradiction.
If $y=a,$ then $t_at_x=t_xt_a=t_at_z$, so that $x=z$, a contradiction.

If $x \cdot y=1$, then by Lemma \ref {lem1} applied to  $t_xt_y=t_at_z$ we have $[x,y] \in \{ [a,z],[z,a+z],[a+z,a]    \}$.
If $x=a,y=z$, then we have a contradiction.
Also $x=z,y=a+z$, gives  a contradiction.
Lastly,  $x=a+z,y=a$, gives  case  (iv). \qed
\medskip

\section {The functions $f_1,f_2,f_3,f_4$}

Given $C \subseteq \mathcal T$ we define four functions $f_i=f_i(C):C \to \mathbb Z^{\ge 0}, i=1,2,3,4$:

\noindent For $a \in C$ we let $f_1(C,a)$ be the number of $b \in C$ such that $a \cdot b=1$ and $a+b \in C$. 

\noindent For $a \in C$ we let $f_2(C,a)$ be the number of $b \in C$ such that $a \cdot b=1$ and $a+b \notin C$. 

\noindent For $a \in C$ we let $f_3(C,a)=|C|$. 

\noindent For $a \in C$ we let $f_4(C,a)$ be the number of $ b \in C, b \ne a,$ such that $a \cdot b=0$. 

For $t_a,t_b \in C$ we will say that {\it the pair $t_a,t_b$ has type $f_1$ if   $a \cdot b=1$ and $a+b \in C$}. Similarly for types $f_2$ and $f_4$. 

Note that for fixed $t_a \in C$ and $t_b \in C, b\ne a,$  either $a \cdot b=1 $ or $a \cdot b=0$. Thus
\begin{align} \label{eq123456}
2f_1(C,t_a)+f_2(C,t_a)+f_4(C,t_a)+1=f_3(C,t_a).
\end{align}

One of our goals will be to show that if $C\subseteq \mathcal T$ is a p-set of a commutative Schur ring $\mathfrak S$ over $G_n$ containing $\mathcal T$, then the functions $f_i(C)$ are constant on $C$.

For $i=1,2,4$ and  a p-set $C$ we let 
\begin{align} \label{eq12345}
(C^2)_{f_i}=\{t_at_b:t_a,t_b \in C \text { where } t_a,t_b \text { has type } f_i\}.
\end{align}

 \begin{lem}\label{lemplus}
For any $C \subseteq \mathcal T$ each $(C^2)_{f_i},i=1,2,4,$ is in $\mathfrak S$.

 \end{lem}
 \noindent {\it Proof}   By Lemma \ref {lem31} 
 $
 \overline{C}^2=|C|+   3(C^2)_{f_1}+ 1(C^2)_{f_2}+2(C^2)_{f_4} $ and we are done.
 \qed 
 \medskip

Now let
$n\ge 4, C \subseteq \mathcal T, |C|>0, C'=\mathcal T \setminus C, \overline C \in \mathfrak S$, so that $\overline {C'} \in \mathfrak S$. 
Let  
 $$ 
 A={\rm Tr}^3(C), \,\, B={\rm Tr}^2(C),\,\, A'={\rm Tr}^3(C'),\,\, B'={\rm Tr}^2(C').
 $$
 Then by Lemma \ref {lem32} 
 $$
 \sum_{T \in A} T^* \in \mathfrak S, \,\,
  \sum_{T \in B} T^* \in \mathfrak S, \,\,
   \sum_{T \in A'} T^*\in \mathfrak S, \,\,
    \sum_{T \in B'} T^* \in \mathfrak S.
 $$
 
 Note that if $U=\{a,b,a+b\} \in A, a\cdot b=1$ so that $U^*=t_at_b+t_bt_a$, then  
 \begin{align}
 \label{eq41} U^*(t_a+t_b+t_{a+b})=2\cdot (t_a+t_b+t_{a+b})
 \end{align}
 while if $U=\{a,b\} \in B,$ where $a,b \in C, a\cdot b=1$, then
  \begin{align}
& \label{eq42} U^*(t_a+t_b)=2 t_{a+b}+t_a+t_{b};\\ 
&U^*(t_a+t_{a+b})=2t_b+t_a+t_{a+b}; \label{eq43}\\
&U^*(t_b+t_{a+b})=2t_a+t_b+t_{a+b}. \label{eq44}
 \end{align}
 
 If $a \cdot b=1$, then we define $\sigma (\{a,b,a+b\}):=t_a+t_b+t_{a+b}$.
 
 \begin{lem} \label{lem42}
 For $A, A'$ as above we have:
 
 \noindent (i) $\sum_{\{a,b,a+b\}\in A} \sigma(\{a,b,a+b\}) \in \mathfrak S;$
 
 \noindent (ii) $\sum_{\{a,b,a+b\}\in A'} \sigma(\{a,b,a+b\}) \in \mathfrak S.$
  \end{lem}
  \noindent{\it Proof} 
  (i) From $\sum_{T \in A} T^*, \overline {(t_1t_n)^G}  \in \mathfrak S$ we get
   $\sum_{T \in A} T^*\cdot  \overline {(t_1t_n)^G}  \in \mathfrak S$. By Eq. (\ref {eq41})  the coefficient of $t_a$ in this product is twice the number of 
   $\{a,b,a+b\}_{b \in V^\sharp}$ that are in $A$. So (i) follows and (ii) is similar, using Eqs.  (\ref {eq42}), (\ref {eq43}), (\ref {eq44})   .\qed
   \medskip

   \begin{cor} \label{corlem42}
   If $C\subset \mathcal T$ is a p-set,  then for   $t_a \in C$ the number of $t_b \in C, a\cdot b=1,$ where $t_{a+b} \in C$,  is independent of $a$,
   and  the function $f_1(C)$ is constant on $C$.
   \end{cor} 
     \noindent{\it Proof}  Let $C$ be a p-set.
     Let $q:=\sum_{\{a,b,a+b\}\in A} \sigma(\{a,b,a+b\}) \in \mathfrak S.$
     The transvection $t_a\in C$ occurs with coefficient $\lambda$ in $q$ if and only if the number of $\{a,b,a+b\}\in A$ with $a$ constant and $b$ varying is $\lambda$.
     Since $C$ is a p-set and $q \in \mathfrak S$, this number must be constant over $C$. But this number is  $f_1(C,t_a)$.\qed
     \medskip

 %  EXAMPLE: For $n=4$ the only subsets satisfying the condition  on a set $C$ that the number of such 2-subspaces contained in $C$ and containing  a given $a$ is independent    of $a$,  are the ones that we obtain using the SGP subgroups. 
   
  We now consider $f_2(C,t_a)$.

\begin{lem}\label{lempr12}
Let $C \subseteq \mathcal T$ and  $t_a \in C$. Then the coefficient of $t_a$ in 
$\overline{C} \cdot (t_{1}t_{2})^G$ is $f_4(C,t_a)$.
\end{lem}
\noindent {\it Proof} This follows directly from Lemma \ref {lem43} as an element of such a product has the form
$t_u\cdot t_vt_w$, where $v \cdot w=0$ and $t_u \in C$.\qed
\medskip

\begin{prop}\label{propf1f2}
Given a p-set $C\subseteq \mathcal T$ the functions $f_1(C), f_2(C), f_3(C), f_4(C)$ are constant.
\end{prop}
\noindent {\it Proof} From Lemma \ref {lempr12}  $f_4(C)$ is  constant  on $C$; also  $f_3(C)=|C|$. 
By  Corollary \ref {corlem42}  $f_1(C)$ is constant on $C$.
But then Eq. (\ref {eq123456})    shows that $f_2(C)$ is constant on $C$.\qed
\medskip

Since  $f_i(C_j,\_\_)$ is constant on each p-set $C_j$ we will write $f_i(C_j)$  for this constant value. 
Define
  \begin{align} \label{eqd1d2}
  &D_1(C):= \sum_{\begin{matrix} T \in B \text { where}\\ T\cap C=\{u,v\}\end{matrix}} t_u+t_v,\quad   
   D_2(C):=  \sum_{\begin{matrix} T \in B \text { where}\\ T\cap C=\{u,v\}\end{matrix}}  t_{u+v}. 
  \end{align}

  \begin{lem} \label{lem44} Let $C\subseteq \mathcal T$ be a  p-set of  $\mathfrak S$. Then   $D_1=D_1(C), D_2=D_2(C) \in \mathfrak S$.
  \end{lem}
    \noindent{\it Proof} 
    From $ \sum_{T \in B} T^*, \overline {C}  \in \mathfrak S $
  we get 
   $\left (  \sum_{T \in B} T^*\right ) \cdot  \overline {C}  \in \mathfrak S.$
   From Eq. (\ref {eq42}) we get 
   $$
   \left (  \sum_{T \in B} T^*\right ) \cdot  \overline {C} =D_1+2D_2 +X,
   $$
   where $X$ does not contain any elements of $\mathcal T$.
 It follows that
   $D_1+2D_2  \in \mathfrak S$.
   
   Now from the definition of $D_1$ we see that $D_1 \subseteq C$, and since $C$ is a p-set  it follows that $C \subseteq D_1$. Thus $D_1=C$ and the result follows.\qed\medskip

\section{$f_2(C)=0$ for some p-set $C$}

Given $\varepsilon \in \{0,1\}=\mathbb F_2$ and any $C \subseteq \mathcal T$ we define the graph $\Gamma^\varepsilon(C)$ to have vertices $C$ and an edge $\{a,b\}\subseteq C, a \ne b,$ if $a \cdot b=\varepsilon$. 
If $a,b \in V^\sharp, a \cdot b=1$, then we will call the set $\{a,b,a+b\}$ a {\it zero-triangle}.

\begin{prop} \label{propuv} Let $C\subseteq \mathcal T$ be a p-set where $f_2(C)=0$ and 
  let  $t_a \in C$. Then $C^{t_a}=C.$    In particular, $\langle C\rangle$ acts by conjugation on  $C$.

  \end{prop}
  \noindent{\it Proof} So assume that $f_2(C)=0$. 
  If we also have $f_1(C)=0$, then  any $t_a \in C$ commutes with every  $t_b \in C$.
  Thus $C^{t_a}=C$ in this case.
  
   Now we assume     $f_1(C)\ne 0$, so that  each $t_{a+b} \in C$ is in a typical  zero-triangle  $\{t_a,t_b,t_{a+b}\}\subseteq C$ (as    $f_1$ is constant on $C$). 
Consider without loss of generality  the conjugation action of $t_{a+b}=t_b^{t_a}=t_a^{t_b}\in C$  
  on three types of vertices of $\Gamma^1(C)$:
  
\noindent (i) those in   $\{t_a,t_b,t_{a+b}\}$;
  
\noindent  (ii) those not in  $\{t_a,t_b,t_{a+b}\}$, but which are   adjacent in $\Gamma^1(C)$ to  any of the vertices in   $\{t_a,t_b,t_{a+b}\}$;
  
\noindent  (iii) those not adjacent in $\Gamma^1(C)$ to any of the vertices   $t_a,t_b,t_{a+b}$. 
  
For (i) the three choices give: 
$
t_a^{t_{a+b}}=t_b \in C, t_b^{t_{a+b}}=t_a \in C, t_{a+b}^{t_{a+b}}=t_{a+b}\in C.
$

For (ii) there are three cases for such a vertex $t_c\in C$: 
\begin{enumerate}

\item $t_c\in C$ is adjacent to $t_a$ in $\Gamma^1(C)$; 

\item  $t_c\in C$ is adjacent to $t_b$ in $\Gamma^1(C)$; 

\item  $t_c\in C$ is adjacent to $t_{a+b}$ in $\Gamma^1(C)$.
\end{enumerate}

\noindent {\bf CASE (ii) 1}: Here $c \cdot a=1$. Since $f_1(C)>0, t_a,t_c \in C,$ and $ f_2(C)=0$ we have
    $\{t_a,t_c,t_{a+c}\}\subseteq C$. Then
     $c \cdot a=1$ shows that either (a) $c \cdot b=1$ or (b) $c \cdot(a+b)=1$. 

If  (a), then $f_2(C)=0$ gives   $\{t_b,t_c,t_{b+c}\}\subseteq C$ is a zero triangle.   Then  $c\cdot (a+b)=0,$ and so 
  $
  t_c^{t_{a+b}}=t_c \in C.
  %,  t_{a+c}^{t_{a+b}}=t_{c+b}\in C.
  $
  
If  (b), then $c\cdot b=0$ and so $f_2(C)=0$ shows that   $\{t_c,t_{a+b},t_{a+b+c}\}\subseteq C$ is a zero triangle.
 Then $c \cdot(a+b)=1$ gives
$
t_c^{t_{a+b}}=t_{a+b+c}\in C.
%,\quad t_{a+b}^{t_{a+b}}=t_{a+b}\in C,\quad t_{a+b+c}^{t_{a+b}}=t_c\in C.
$

\noindent {\bf CASE (ii) 2}: This case is  analogous to CASE 1 (just switch $a$ and $b$).

\medskip 

\noindent {\bf CASE (ii) 3}: Here  $c \cdot (a+b)=1$. Then, since $t_c,t_{a+b} \in C, c \cdot (a+b)=1$, 
 again we have
$
t_c^{t_{a+b}}=t_{a+b+c} \in C.
$ 

For (iii): here  $c\cdot a=c\cdot b=0$ and so 
$t_c^{t_{a+b}}=t_c$ . Thus in all cases we have $t_c^{t_{a+b}}\in C$.
This gives the result.
\qed
\medskip

%\begin{lem}\label{corsym} Let $\{t_a,t_b\}$ be an edge of $\Gamma^1(C)$ so that $a \cdot b=1$. Then there is $\alpha \in \langle C\rangle$ such that $t_a^\alpha=t_b, t_b^\alpha=t_a$ and $\alpha^2=1$.

%In particular, if $\Gamma^1(C)$ is connected, then  $ \langle C\rangle$ acts transitively on the vertices of $\Gamma^1(C)$.\end{lem} 
%\noindent  {\it Proof} We see that in this situation $\{t_a,t_b,t_{a+b}\}$ is a zero-triangle in $C$. Then take $\alpha=t_{a+b}=t_at_bt_a\in \langle t_a,t_b\rangle\le \langle C\rangle$.\qed\medskip

\begin{lem}\label{lemsizes}
(i) The class of $t_{1}t_{2}$ has size $(2^n-1)(2^{n-2}-1)$.

\noindent (ii) The class of $t_{1}t_{n}$ has size $(2^n-1)2^{n-1}/3$.

\noindent (iii) The number of dimension $2$ symplectic  subspaces of $V$ is equal to the number of zero triangles $\{a,b,a+b\}$ in $V$, which is equal to 
$(2^n-1)2^{n-2}/3$.

\noindent (iv)  The order of $G_n, n=2m,$ is $2^{m^2}(2^n-1)(2^{n-2}-1)(2^{n-4}-1)\cdots (2^2-1).$
\end{lem} 
\noindent{\it Proof} (i) This is equal to the number of subsets $\{a,c\} \subset V^\sharp$ of size $2$  where $a\cdot c=0$. The rest are similar.\qed

\section {An overview of the proof}

Given any subsets $C_1,\cdots,C_r\subseteq \mathcal T$   we define the graph $ \Gamma^1(C_1,\cdots,C_r)$ to have vertices $C_1,\cdots,C_r$ and an edge between $C_i, C_j, i \ne j,$ if there are $t_a \in C_i, t_b \in C_j$ that do not commute. 
If  a partition  $\mathcal T=\cup _{i=1}^r C_i$ is understood (as given by  the p-sets of a Schur ring for example), then we will denote 
 $ \Gamma^1(C_1,\cdots,C_r)$ by  $ \Gamma^1(\mathcal T)$.
 The following are steps in the proof:
 
 \begin{enumerate}
 
\item  [Step 1:] Show that  $ \Gamma^1(\mathcal T)$ is connected and bipartite: Corollary \ref {corK}.

\item  [Step 2:]  Show that we cannot have an edge $\{C_i,C_j\}$ of $ \Gamma^1(\mathcal T)$ where $f_2(C_i), f_2(C_j)>0$: Theorem \ref {thmf1g1222}.
 
\item   [Step 3:]  Show that we can assign an orientation to each edge of $ \Gamma^1(\mathcal T)$, so that 
 $ \Gamma^1(\mathcal T)$ is a directed graph with an edge $C_i \to C_j$ if $f_2(C_i)>0, f_2(C_j)=0$: Lemma \ref {lemorientgr}.
 
  \item   [Step 4:]  If we have $i \ne j$ with either $f_2(C_i),f_2(C_j) >0$ or  $f_2(C_i)=f_2(C_j) =0$, then $C_i \perp C_j:$ Corollary \ref {corgr111} and Proposition \ref {propf2f2c}.
 
 \item   [Step 5:]  We show that we cannot have distinct vertices $C_i, C_j, C_k$ of  $ \Gamma^1(\mathcal T)$ with directed edges $C_i \to C_j$ and $C_k \to C_j$ of
  $ \Gamma^1(\mathcal T)$: Proposition \ref {lem222}.

 \item   [Step 6:] The graph $ \Gamma^1(\mathcal T)$ is a connected star graph with one source vertex.
   More specifically this gives  what we call the {\it standard situation}: the p-sets $C_k$ can be numbered so that $C_1$ is the source and  so that $f_2(C_1)\ge1, f_2(C_i)=0,1<i\le r$. Further we also have:  $C_i\perp C_j$ for distinct $i,j >1$ and there is a directed edge $C_1\rightarrow C_i$ for each $1<i\le r$: Corollary \ref {cor 222}. 
   
  \item   [Step 7:]   Each $C_i, 1\le i\le r,$ is invariant under the action by conjugation of the group $H:=\langle C_2,\cdots,C_r\rangle$:
  Corollary \ref {lemc1c2comm}.
      
 \item   [Step 8:]   No more than one of the $C_i, i>1,$ has $f_1(C_i)>0$: Proposition \ref {prop3cpts}.

 \item  [Step 9:] For $n\ge 6, r>1,$ the case where each $C_i,1<i\le r,$ is isotropic does not occur:
 Propositions \ref {propisotropic1}, \ref {propisotropic2}, \ref {propisotropic3}.
 
 By Step 9 we can assume that $f_1(C_2)>0$ and $C_3,\cdots,C_r$ are isotropic (the {\it up-dated standard situation}).

 %  Let $X=\cup_{i=2}^r C_i$; then $X={\rm Span}(X)^\sharp$. NEEDED and ${\rm Span}(X)=V$

\item  [Step 10:] For $n\ge 6$ the standard situation    is only possible if either (A) $r=2$; or (B) $r=3$ and $|C_3|=1$: Corollary \ref {cor32111}.

\item  [Step 11:]  The graph $\Gamma^1(C_2)$ is connected and  $H=\langle C_2\rangle$ is an irreducible $3$-transposition subgroup of $G$: Propositions \ref     {propconn} and \ref {lemconnC2}.

\item  [Step 12:]   Here we obtain relations among the `variables' $a_i=f_i(C_1), b_i=f_i(C_2), 1\le i\le 4,\l_1,\l_2$: see $\S 13$. 
%\medskip

Now the results of McLaughlin \cite {McL} show that any  proper   irreducible subgroup of ${\rm Sp}(n,2)$ generated by transvections is one of:
$$  (i) \,\, S_{n+1};\quad 
 (ii)  \,\, S_{n+2};\quad  (iii)  \,\, {\rm SO}^+(n,2);
\quad  (iv)  \,\,  {\rm SO}^-(n,2).
$$

% \noindent (v) ${\rm Sp}(n,2)$.

Here irreducible means that there are no invariant subspaces. This gives the possibilities for $H.$
The next steps are to  consider  these four cases separately.

%\medskip 

\item  [Step 13:]  If we have (A), and  $H$ is a symmetric group (cases (i) and (ii) above), then we obtain a contradiction from   the set of relations  that we have obtained (if $n\ge 14$).
 
 We also similarly show that (B) does not happen for any  of (i)-(iv).

 \item  [Step 14:] We then show what happens when we have (A) and $H={\rm SO}^\pm(n,2)$.  

     \end{enumerate}

 \section{ Step 1: $ \Gamma^1(\mathcal T)$ is connected}

Let $\mathcal T=\bigcup_{i=1}^r C_i$, a union of p-sets.
  We first   show that   $ \Gamma^1(\mathcal T)$ is connected:

\begin{lem}\label{lemK}
Define $v_1,v_2,\cdots,v_{n+1}\in V^\sharp$ as follows:
\begin{align*}
&v_1=e_1, \,\,\, v_2=e_n\\&
v_3=(e_1+e_2)+(e_n),  \,\,\, v_4=(e_1)+(e_{n-1}+e_n),\\&
v_5=(e_1+e_2+e_3)+(e_{n-1}+e_n),  \,\,\, v_6=(e_1+e_2)+(e_{n-2}+e_{n-1}+e_n),\\&
v_7=(e_1+e_2+e_3+e_4)+(e_{n-2}+e_{n-1}+e_n), \\&  v_8=(e_1+e_2+e_3)+(e_{n-3}+e_{n-2}+e_{n-1}+e_n),\\&
\qquad \vdots\\&
v_{n+1}=e_1+e_2+\cdots+e_n.
\end{align*}
Then $v_i\cdot v_j=1$ for all $1\le i\ne j\le n+1$. Further $\{v_1,\cdots,v_{n}\}$ is a basis.
\end{lem}
\noindent {\it Proof} Just check that $v_i\cdot v_j=1$ for $1\le i\ne j\le n+1$.
\qed
\medskip

\begin{cor} \label{corK}
 The graph $\Gamma^1(\mathcal T)$ is connected.
\end{cor}
\noindent {\it Proof} Let $v_1,\cdots,v_{n+1}$ be as in Lemma \ref {lemK}. Then   $\Gamma^1(\{v_1\},\{v_2\},\cdots,\{v_{n+1}\})$   is a complete graph and so is connected.
 Since $\{v_1,\cdots,v_{n}\}$ is a basis for $V$, if $v \in V, v \ne 0$, then there is some $1\le i\le n$ such that $v_i \cdot v=1$. Thus 
    $\Gamma^1(\mathcal T)$ is connected.\qed
    \medskip

\begin{prop} \label{propf2f233} Let $\mathcal T=\bigcup_{i=1}^r C_i$ as a union of p-sets.
Then we cannot  have  $f_2(C_i)=0$ for all $1\le i\le r$.
\end{prop}
\noindent {\it Proof}  If $f_2(C_i)=0$ for all $1\le i\le r$, then the graph $\Gamma^1(\mathcal T)$ is  not connected.\qed
\medskip

  \section {Steps 2, 3, 4: $f_2(C_i)f_2(C_j)\ne0$ cannot happen}

 \begin{thm}\label{thmf1g1222}
If $\mathfrak S$ is a commutative Schur ring, then we cannot have an edge $\{C_i,C_j\}$ of $\Gamma^1(\mathcal T)$  where $f_2(C_i), f_2(C_j) >0$.
    \end{thm}
    \noindent {\it Proof}   So assume that  we have distinct p-sets $C_1,C_2$ giving an edge of $\Gamma^1(\mathcal T)$ where $f_2(C_1), f_2(C_2) >0,$ so that 
     there are $t_a,t_b \in C_1, a\cdot b=1,$  with $t_{a+b} \in C_2$. 
     %Here the pair $t_a, t_b\in C_1$ is of type $f_2$.
    
   Since $C_1, C_2$ are p-sets,  $f_2$ is constant on $C_1$ and  on $C_2$.     
     Since $t_{a+b}\in C_2$ and $f_2(C_2) \ne 0$ there is $c \in C_2$  such that $t_{a+b},t_c\in C_2$ has type $f_2$, so that $c \cdot (a+b)=1$ and $t_{a+b+c} \notin C_2$.
      There are two cases:  either
     
\noindent     {\bf (A):} $a+b+c \in C_1$; or 
     
\noindent    {\bf  (B):} there is $m>2$ such that  $a+b+c \in C_m $. We may assume that $m=3$.
    
    In either of these cases we have the following: 
If $a\cdot (a+b+c)=0=b \cdot(a+b+c)$, then $a\cdot c=1=b \cdot c$, showing $c\cdot (a+b)=0$, a contradiction. So either $a\cdot (a+b+c)=1$ or 
$b \cdot (a+b+c)=1$. Since we can interchange $a$ and $b$ in what we have done so far, we may, without loss of generality, assume that 
$b \cdot (a+b+c)=1$.
Then:
$$
a\cdot b=1,\,\,\,  b \cdot c=0, \,\,\, a\cdot c=1.
$$

Since $b \cdot (a+b+c)=1$  the pair $b,a+b+c$ is either type $f_1$ or type $f_2$. There are now four  cases to be considered:

\noindent {\bf Case A1:}  (A) and $b,a+b+c$ has type $f_1$. Here $\{ a,b,a+b+c,a+c\} \subseteq C_1, \{a+b,c \} \subseteq C_2$.

\noindent {\bf Case A2:}   (A) and   $b,a+b+c$ has type $f_2$. Here $(a+b+c)+b=a+c \notin C_1$ and so there are two cases that we can index  as follows:
  
\noindent {\bf Case A21:}  $a+c \in C_2$, so that
 $\{a,b,a+b+c\}\subseteq C_1, \{a+b,c,a+c \} \subseteq C_2$. 

\noindent {\bf Case A22:}  $a+b \in C_3$, so that
 $\{a,b,a+b+c\}\subseteq C_1, \{a+b,c \} \subseteq C_2, a+c \in C_3$. 
 Here  $a \in C_1, c \in C_2, a+c \in C_3, a\cdot c=1,$ which violates Lemma \ref{lem3}, so this case does not happen.

\noindent {\bf Case B1:}  (B)  and the pair $b,a+b+c$ has type $f_1$. This case cannot happen since $b \in C_1, a+b+c\in C_3 \ne C_1$.

\noindent {\bf Case B2:}   (B) and the pair $b\in C_1,a+b+c\in C_3$ has type $f_2$.  Here there are two cases (by Lemma \ref {lem3}):

\noindent  {\bf Case B21:} we have (B2) and $a+c \in C_1$;

\noindent {\bf   Case B22:} we have (B2) and $a+c \in C_3$.  
Here  $t_a \in C_1, t_c \in C_2, t_{a+c}\in C_3$ and $a \cdot c=1$,   violating Lemma \ref {lem3}.

\medskip

We thus need to consider the  three  cases A1, A21, B21, as follows:

First,  from $
a\cdot b=1,  b \cdot c=0,  a\cdot c=1
$   we can assume: $a=e_1,b=e_n,c=e_2+e_n$. Let $W={\rm Span}(e_1,e_2,e_n)$. Note that 
$$
a \cdot (b+c)=b\cdot (b+c)=c \cdot (b+c)=0,
$$
 so that $t_{b+c}$ commutes with each of
$t_a,t_b,t_c,t_{a+b},t_{a+c}, t_{a+b+c}$.

\noindent {\bf Case A1} Here 
\begin{align}
\label{eqnA1}&
\{ a,b,a+b+c,a+c\} =\{e_1,e_n,e_1+e_2,e_1+e_2+e_n\}
\subseteq C_1,\\& \{a+b,c \}=\{e_1+e_n,e_2+e_n\}  \subset C_2.\notag 
\end{align}

Lemma \ref {lem32} gives $ (C_1^2)_{f_2} \in \mathfrak S$.
 Now $t_{2,n} \in C_2, t_1t_n \in  (C_1^2)_{f_2}$  and so 
$t_{2,n} t_1t_n\in C_2\cdot  (C_1^2)_{f_2}\subseteq C_2C_1^2$; we will  show that
$t_{2,n} t_1t_n\notin  (C_1^2)_{f_2}C_2\subseteq C_1^2\cdot C_2$.
 We need:

\begin{lem} \label{lemA1case111}
If $t_u,t_v,t_w \in \mathcal T$,  $t_ut_vt_w=\a:=t_{2,n}t_1t_n$, then $(t_u,t_v,t_w)$ is one of:

\noindent (i) $(t_{1,2},t_{1},t_{1,2,n})$   \,\,\, (ii) $(t_{1,2},t_{1,2,n},t_{2,n})$; \,\,\, (iii) $(t_{1,2},t_{2,n},t_{1})$;

\noindent (iv) $(t_{1},t_{1,2},t_{1,2,n})$, \,\,\, (v) $(t_{1},t_{n},t_{1,2})$;\,\,\, (vi) $(t_{1},t_{1,2,n},t_{n});$

\noindent (vii) $(t_{n},t_{1,2},t_{2,n})$; \,\,\, (viii) $(t_{n},t_{1,n},t_{1,2})$;\,\,\,(ix) $( t_{n},t_{2,n},t_{1,n})$;

\noindent (x) $(t_{1,2,n},t_{n},t_{2,n})$; \,\,\, (xi) $(t_{1,2,n},t_{2,n},t_{n})$;

\noindent (xii) $(t_{1,n},t_{1,2},t_{1});$ \,\,\,(xiii) $(t_{1,n},t_{1},t_{1,2})$;

\noindent (xiv) $(t_{2,n},t_{1},t_{n})$;\,\,\, (xv) $(t_{2,n},t_{n},t_{1,n})$;\,\,\, (xvi) $(t_{2,n},t_{1,n},t_{1})$.

In particular $u,v,w \in \{t_{1},t_{n},t_{1,n},t_{1,2},t_{2,n},t_{1,2,n}\}$. 

\end{lem}
\noindent {\it Proof} First one checks that each of these cases gives a solution to  $t_ut_vt_w=\a$.
Next we note that $(e_n)\a=e_n+(e_1+e_n),
(e_1)\a=e_1+(e_1+e_2+e_n)$ and $(e_{n-1})\a=e_{n-1}+(e_1+e_2)$, showing that 
${\rm Span}(e_1,e_2,e_n)= {\rm Span}(e_1+e_n,e_1+e_2+e_n,e_1+e_2)\subseteq {\rm Span}(u,v,w)$. Thus {\rm Span}$(u,v,w)={\rm Span}(e_1,e_2,e_n)$.
Thus $u,v,w$ is a basis for  {\rm Span}$(e_1,e_2,e_n)$, and we check all such sets to find the $16$ cases listed.\qed
\medskip

Now  in case A1  we have $t_{2+n} t_1t_n\in C_2\cdot  (C_1^2)_{f_2}$. Since $\mathfrak S$ is commutative 
 there are $t_u, t_v \in C_1$ of type $f_2$ and $t_w \in C_2$ such that
$t_ut_vt_w=t_{2,n}t_1t_n$. By Lemma \ref {lemA1case111} 
$$
t_u,t_v,t_w \in S:=\{t_1,t_n,t_{1,n}, t_{2,n},t_{1,2,n},t_{1,2}\}.
$$
In fact from Eq. (\ref {eqnA1})  we know  which $C_i$ each of the elements of $S$ is contained in.
Using  Lemma \ref {lemA1case111} we now show that there are no  triples $t_u,t_v \in C_1$ of type $f_2$ and $t_w\in C_2$ satisfying these conditions. 
Specifically, in Lemma \ref {lemA1case111}:  

\noindent options (i) and (iv)  are not possible, since $t_w=t_{1,2,n} \notin C_2$; 

\noindent options (ii), (vii), (x)   and (xiii) are not possible, since $t_u,t_v$   is type $f_1$; 

\noindent options (iii), (xii) and (xvi)  are not possible, since $t_w=t_{1} \notin C_2$; 

\noindent options (v) and (viii) are  not possible, since $t_w=t_{1,2} \notin C_2$; 

\noindent options (vi), (xi) and (xiv) are  not possible, since $t_w=t_{n} \notin C_2$; 

\noindent options (ix) and (xv) are not possible since $t_u\cdot t_v=0$.

Thus Case A1 does not happen.
\medskip

\noindent {\bf Case A21}: Here 
\begin{align}
& \label{lemA2}
\{a,b,a+b+c\} =\{e_1,e_n,e_1+e_2\} \subseteq C_1,\\&
 \{a+b,c,a+c\}=\{e_1+e_n,e_2+e_n,e_1+e_2+e_n\}   \subseteq C_2,\notag 
 \end{align}
and we have  
$
 (C_1^2)_{f_2}\cap \mathcal T(W)^2=\{t_1t_n,t_nt_{12n}\}.
$

Now we consider $\b:=t_{1,n}\cdot t_{n}t_{1,2} \in C_2 \cdot [ (C_1^2)_{f_2}\cap \mathcal T(W)]$ where commutativity means that there are  $t_u,t_v \in C_1$ of type $f_2$ and $t_w \in C_2$ such that $t_ut_vt_w=\b.$ 
Now we show that $\b \notin [ (C_1^2)_{f_2}\cap \mathcal T(W)]\cdot C_2$ 
using  a result similar to   Lemma  \ref {lemA1case111} to  conclude  that there are no  $u,v,w$ satisfying the required conditions: $ t_u,t_v \in C_1$ of type $f_2$,
$ t_w \in C_2$ and  $t_ut_vt_w=\b.$  This shows that Case A21 does not occur. 
\medskip

\noindent {\bf Case B21}: 
Here $\{a,b,a+c\}=\{t_1,t_n,t_{1,2,n} \} \subseteq C_1, \{a+b,c\}=\{t_{1,n},t_{2,n}\}\subseteq C_2, \{a+b+c\}=\{t_{1,2}\} \subseteq C_3$. 

Now let $U:=C_1C_2 \cap (t_1t_2)^G \cap {\rm Sp}(W) = \{t_nt_{2+n}, t_{1+2+n}t_{1+n}\}$.
 We will need
 
 \begin{lem}\label{lemcaseb21}
 Let $\a=t_{n}t_{2,n}t_{1,n}$. Then the only $t_u,t_v,t_w \in \mathcal T$ such that $v \cdot w=0$ and $t_ut_vt_w=\a$ are
 
 \noindent (i) $t_{1,2,n}, t_{n}, t_{2,n}$;

 \noindent (ii) $t_{1,n}, t_{1}, t_{1,2}$;
 
  \noindent (iii) $t_{1,2,n},  t_{2,n}, t_{n}$;

  \noindent (iv) $t_{1,n}, t_{1,2}, t_{1}.$
 \end{lem}
 
 \noindent {\it Proof} This is similar to the proof of Lemma \ref {lemA1case111}. First note that
$$ e_1\a=e_1+(e_1+e_2+e_n),\quad 
e_n \a=e_n +(e_1+e_n),\quad e_{n-1}\a=e_{n-1}+(e_1+e_2).
$$
Thus {\rm Span}$(e_1+e_2+e_n, e_1+e_n,e_1+e_2\} \subseteq {\rm Span}(u,v,w)$, which shows that $u,v,w$ is a basis for 
{\rm Span}$(\{e_1+e_2+e_n, e_1+e_n,e_1+e_2\}) ={\rm Span}(\{e_1,e_2,e_n\})$. 
This again reduces the proof to a short  enumeration.
\qed\medskip

 Let $\a:=t_{n}t_{2,n}t_{1,n}$ and note that   $t_{n}t_{2,n} \in C_1C_2 \cap (t_1t_2)^G$, where $ E:=C_1C_2 \cap (t_1t_2)^G$ is an $\mathfrak S$-set.
We show that $C_2$ does not commute with $E$; in particular $t_{1,n} \in C_2$ and so
$( t_{n}t_{2,n})\cdot t_{1,n} \in EC_2$. 
We show that $( t_{n}t_{2,n})\cdot t_{1,n} \notin C_2E$.

Since $E$ and $C_2$ commute,  there are 
$t_u \in C_2$ and  $t_vt_w \in E$, so that $v \cdot w=0$ and $t_u(t_vt_w)=\a$.  
But the possibilities for such $u,v,w$ are given in Lemma \ref {lemcaseb21} (i)-(iv).
Now we cannot have (i) or (iii) since we would then have $t_u=t_{1,2,n}\in C_1$, contradicting $t_u \in C_2$.
 We cannot have (ii) or (iv) since $t_{1,2}\in C_3$. 
This concludes the proof of  Case B21 and so completes the proof of Theorem \ref {thmf1g1222}.\qed
\medskip

\begin{prop} \label{propf2f2c}
If $C_i,C_j,i \ne j,$ are commuting  p-sets with $f_2(C_i)=f_2(C_j)=0$, then    $C_i\perp C_j$.
\end{prop}\
\noindent {\it Proof} Suppose that $C_i,C_j,i \ne j,$ are   p-sets with $f_2(C_i)=f_2(C_j)=0$
and  we have $u \in C_i,v \in C_j$ with  $u \cdot v=1$. Since $C_i, C_i$ are p-sets that commute, we must have $u+v \in C_i \cup C_j$ by Lemma \ref {lem3}. So suppose without loss of generality that
$u,u+v \in C_i, v \in C_j$; then this
 contradicts $f_2(C_i)=0$.\qed
\medskip

      We will show how to make  $\Gamma^1(C_1,\cdots,C_r)$  into a directed graph i.e. we need to assign an orientation to each edge.
           So suppose that there is an edge $\{C_i, C_j\}$.
     Then we cannot have $f_2(C_i)=f_2(C_j)=0$ 
     as then $C_i \perp C_j$ by Proposition \ref{propf2f2c}.

     So assume without loss of generality  that $f_2(C_i)>0.$ Then Theorem \ref {thmf1g1222}  shows that $f_2(C_j)=0$ and we orient the edge $\{C_i,C_j\}$ from $C_i $ to $C_j$.  
   Thus 
   \begin{lem} \label{lemorientgr} As defined above, 
    $\Gamma^1(C_1,\cdots,C_r)$ is a directed graph.
    \qed
    \end{lem}

       \begin {lem} \label{lem12bip} The graph $\Gamma^1(C_1,\cdots,C_r)$, considered as an un-directed graph,  is bipartite and connected.
      \end{lem}
         \noindent {\it Proof}  Whenever we have an edge $\{C_i,C_j\}$ exactly one of $C_i, C_j$ has non-zero $f_2$. This determines the bipartite partition of 
      $\Gamma^1(C_1,\cdots,C_r)$. We have already noted that it is connected. \qed
      \medskip

   If $\{C_i, C_j\}$ is an edge of $\Gamma^1(\mathcal T)$ where the orientation is from $C_i$ to $C_j$, then we notate this as
   $C_i \to C_j$.

\begin{cor}\label{corgr111}
If $C_i, C_j,i \ne j,$ are   p-sets with $f_2(C_i), f_2(C_j) >0,$  then  $C_i\perp C_j$.
\end{cor} 
\noindent {\it Proof}
The situation where  $f_2(C_i), f_2(C_j) >0$ and there are $t_a \in C_i, t_b\in C_j$ with $a \cdot b=1$ is prohibited by Theorem \ref {thmf1g1222}.
\qed\medskip

This concludes the proofs of Steps 2, 3, 4.

\section {Steps 5, 6, 7: Star graph}

                   \begin{prop}\label{lem222} The situation where there are distinct $i,j,k$ such that
             $C_i \rightarrow C_j$ and $C_k \rightarrow C_j$ are edges of $\Gamma^1(C_1,\cdots,C_r)$ does not occur.  
             \end{prop}
             \noindent {\it Proof}  
            Since  we have  $C_i \rightarrow C_j$ there are $a \in C_i,a+b \in C_j$ where $a \cdot (a+b)=1$. Since $f_2(C_j)=0$ we see that $b=a+(a+b) \in C_i$ and the pair $a,b\in C_i$ has type $f_2$. 
            Now define $D_2(C_i)\subset \mathcal T$ as in Eq. (\ref {eqd1d2}). Then we have $t_{a+b} \in D_2(C_i) \cap C_j \ne \emptyset$.
            Applying  Lemma \ref {lem44}  we see that $C_j \subseteq D_2(C_i)$ since $C_j$ is a p-set. We similarly obtain
            $C_j \subseteq D_2(C_k)$. 
           So  there are  $ t_c, t_d \in C_k$, of type $f_2$, with
             $c+d=a+b$. Now 
             $$
             c \cdot (a+b)=c \cdot(c+d)=1,\text { and so } 1 \in \{c \cdot a, c \cdot b\};
             $$
             but $c \in C_k, a,b \in C_i$, and so there is an edge $\{C_i, C_k\}$ in $\Gamma^1(\mathcal T)$, which contradicts the fact that $\Gamma^1(C_1,\cdots,C_r)$ is bipartite.\qed
             \medskip

                Let $\Gamma^1(C_1,\cdots,C_r)^+$ be the set of vertices $C_i$ of  $\Gamma^1(C_1,\cdots,C_r)$ where $f_2(C_i)>0$, and let $\Gamma^1(C_1,\cdots,C_r)^0$ denote the rest of the vertices. A {\it directed star graph with one source} is a directed graph $\Gamma$  where there is a vertex $v$ (the source) such that for all vertices $w\ne v$ we have an edge from $v$ to $w$, these being the only edges of $\Gamma$.

             \begin {cor} \label {cor 222} The set  $\Gamma^1(C_1,\cdots,C_r)^+$ has exactly  one element. 
             In particular,  $\Gamma^1(C_1,\cdots,C_r)$ is a directed star graph with one source. 
             
             \end{cor}
                    \noindent {\it Proof} Recall that $\Gamma^1(C_1,\cdots,C_r)$ is connected. Thus  if $\Gamma^1(C)_1,\cdots,C_r)^+$ has two  vertices $C_i, C_k, i \ne k$, then there
                    is a path between them in $\Gamma^1(C_1,\cdots,C_r)$. That path in $\Gamma^1(C_1,\cdots,C_r)$  has to have vertices $C_u, C_v$ in $\Gamma^1(C_1,\cdots,C_r)^+$ at distance $2$. This then gives us the situation described in Lemma \ref {lem222}, and so we have a  contradiction.  The rest follows. 
              \qed\medskip
              
                        So we now have the standard situation defined in Step 6.

    Now let $t_u \in C_k,1<k\le r$. Then, by Proposition \ref {propuv},  
    $C_k^{t_u}=C_k$.
    We also have $t_u\in C_k$ and $C_k\perp C_j,1<j\ne k\le r$
which  then gives $C_j ^{t_u}=C_j$ for all $1<j\ne k\le r$.
 But $\mathcal T^{t_u}=\mathcal T, \mathcal T=C_1\cup C_2\cup \cdots \cup C_r$, so that $C_1^{t_u}=C_1$. Thus: 
    
    \begin{lem} \label{lemc1c2comm}
     Let $\mathcal T=C_1\cup C_2\cup \cdots \cup C_r$ be in the standard situation and let $t_u \in C_i, 1<i\le r$.
Then     $C_i^{t_u}=C_i$ for all $1\le i\le r$.

In particular, if $H=\langle C_2\rangle \times \langle C_3\rangle \times \cdots \times \langle C_r\rangle\le G,$ then for $h \in H$ we have
$C_i^h=C_i$ for $1\le i\le r$, and the $C_i$ are unions of orbits of the action of $H$ on $\mathcal T$. 
\qed
\end{lem}

\section {Step 8: $f_1(C_i), i>1$}

For $a,b \in \mathcal T$ we will say that {\it $a$ meets $b$} if $a \cdot b=1$. 
       
    \begin{prop} \label{prop3cpts} Assume  the standard situation.
    Then no more than one of the $C_i, i>1,$ can have $f_1(C_i)>0$. 
    \end{prop}
         \noindent {\it Proof} So assume that $f_{1}(C_2), f_{1}(C_3)>0$. 
          Without loss we can assume that $e_1,e_n \in C_2$, which then gives $\{e_1,e_n,e_1+e_n\}\subseteq C_2$ as $f_1(C_2)>0, f_2(C_2)=0$. Since $ f_1(C_3)>0$ and $C_2 \perp C_3$ we may similarly assume that
         $e_2,e_{n-1}, e_2+e_{n-1}\in C_3$. Let $W:={\rm Span}(e_1,e_n,e_2,e_{n-1})$. 
         Now if $a \in C_2, b \in C_3$, then $a$ meets elements of $C_2$ and $b$ meets elements of $C_3$, so that $a+b$ meets elements of $C_2$ and of $C_3$. 
         Since $C_i \perp C_j$ for all $i\ne j\ge 2$ we must   have $a+b \in C_1$. It follows that
 \begin{align*}&
         \{e_1+e_2,e_1+e_{n-1},e_1+e_2+e_{n-1},
         e_n+e_2,e_n+e_{n-1},e_n+e_2+e_{n-1},\\&\qquad\qquad 
         e_1+e_n+e_2,e_1+e_n+e_{n-1},e_1+e_n+e_2+e_{n-1}\} \subseteq C_1.
\end{align*}
Note that we now have: $(C_1 \cap W) \cup (C_2 \cap W) \cup (C_3\cap W)=W^\sharp$ (a disjoint union).

Define the multisets
$$
X=C_1C_2 \cap (t_1t_2)^G, \quad Y=C_1^2 \cap (t_1t_n)^G,\quad Z=C_1\cdot Y.
$$
    First note that $X,Y, Z \in \mathfrak S$ and that 
    $$
    \a_1:=t_{1,2,n-1,n} t_{1,n}\in X, \quad t_{n-1,n}t_{2,n} \in Y,
    $$
     so that
    $\a_2:=t_{1,n-1,n}t_{n-1,n}t_{2,n} \in Z$.   
    Let $\a:=\a_1\a_2\in XZ$.
    
    \begin{lem}\label{lemalpha} Assume that $n\ge 6$.  Then
    
\noindent   (i) As a product of elements of $\mathcal T$ $\alpha$ has minimal  length $5$.
    
 \noindent   (ii)     If we write $\a=t_at_b t_ut_vt_w$, then $a,b,u,v,w \in W$.
    
  \noindent  (iii) $\a$ occurs in $XZ$ with multiplicity $5$.
    
  \noindent      (iv) $\a$ occurs in $ZX$ with multiplicity $3$.

    \end{lem}
           \noindent {\it Proof}    Note that $\a=\a_1\a_2 \in XZ \subset  C_1C_2 \cdot C_1^3$ and so $\a$ has minimal length at most $5$ as a product  of elements of $\mathcal T$. 
    
 \noindent (i)    Write $\a=t_at_b t_ut_vt_w$, where we allow some of the $a,b,u,v,w$ to be zero, with  the convention  that $t_0=id$. Now note that $\{v\a-v:v \in W\}=W$ and so $W \subseteq {\rm Span}(a,b,u,v,w)$, 
 where $\dim  {\rm Span}(a,b,u,v,w)\le 5$. If   $\dim  {\rm Span}(a,b,u,v,w)=4$, then
     ${\rm Span}(a,b,u,v,w)=W$, and $a,b,u,v,w \in W$.
      If $\dim    {\rm Span}(a,b,u,v,w)  =5$,  then
    ${\rm Span}(a,b,u,v,w)={\rm Span}(e_1,e_2,e_{n-1},e_n,c)$ for some $c \notin W$; further we can change $c$ to $c+u, u \in W$ so that $c\perp W$.     By a change of basis we can then assume that $c=e_3$ (since $n\ge 6$). 
    There are $31$ transvections in $\mathcal T({\rm Span}(e_1,e_2,e_{n-1},e_n,e_3))$ and so we can just find all  $k$-tuples, $1\le k\le 5$, of such transvections whose product is $\a$. Doing this shows that we only have $k=5$, giving (i), and also showing  that  $a,b,u,v,w \in W$, giving  (ii).

    \noindent (iii)  The above has reduced the situation to where we are able to write $\a=t_at_b t_ut_vt_w$, with $t_a,t_b, t_u,t_v,t_w \in \mathcal T(W)$. 
    If we choose $\b_1=t_at_b \in X, \b_2=t_ut_vt_w \in Z$, then another short calculation shows that there are $3$ such  pairs 
    $[\b_1,\b_2]$ with $\b_1=t_at_b \in X \cap \mathcal T(W)^2, \b_2=t_ut_vt_w \in Z \cap \mathcal T(W)^3$ and  $\b_1\b_2=\a$. 
 
   For   pairs $[\b_1,\b_2]$ with $\b_1=t_ut_vt_w \in Z \cap \mathcal T(W)^3, \b_2=t_at_b \in X \cap \mathcal T(W)^2$ and  $\b_1\b_2=\a$ we have
   $5$ such products. This gives (iii) and (iv).
What we have shown is that
   $X,Z \in \mathfrak S$, but $XZ \ne ZX$ as multisets;  Lemma \ref {lemalpha} follows.\qed
   \medskip
   
The contradiction to commutativity obtained from   Lemma \ref {lemalpha}  concludes the proof of Proposition \ref {prop3cpts}.\qed
\medskip  

By Proposition \ref {prop3cpts} there is at most one of the $C_i, 1<i\le r,$ that is not isotropic. We   now reorder the $C_i$ so that if there is one  non-isotropic $C_i, 1<i\le r$, then it is $C_2$. We now  consider this condition as being a part of the standard situation.

\section   {Step 9: the Isotropic cases}
There are three cases for Step $10$. We first prove:

\begin{prop} \label {propisotropic1}
  Suppose that we have the standard situation where $C_2,\cdots, C_r,$ $ r\ge 2,$ are all isotropic. Let $X:={\rm Span}(C_2,\cdots,C_r)^\sharp$, and assume that
 $X \ne C_2\cup\cdots\cup C_r$. Then $C_1,\cdots,C_r$ are not the p-sets of a commutative Schur ring.
\end{prop}
 \noindent {\it Proof}
  So assume that  $C_1,\cdots,C_r$ are   the p-sets of a commutative Schur ring, that $C_2,\cdots, C_r$ are all isotropic and let $X={\rm Span}(C_2,\cdots,C_r)^\sharp$. Let
 $a \in X \setminus  (C_2\cup\cdots\cup C_r)$. Then $a \in C_1$.
But $\Gamma^1(C_1,\cdots ,C_r)$ is connected (as an undirected graph), so   $f_2(C_1)>0$.
Since $a \in C_1$ there is $b \in C_1$ where $a, b$ has type $f_2$, so that $t_{a+b} \in C_2\cup\cdots\cup C_r$ and $a\cdot b=1$.
But   $a \in X, a+b \in X$ and  $X$   isotropic, so   $a\cdot b=1$  is a contradiction.\qed
\medskip

Thus we assume $X={\rm Span}(C_2,\cdots,C_r)^\sharp=\bigcup_{i=2}^r C_i,r\ge 2,$ so that $\dim X \le n/2$ as $X$ is isotropic. The next case is:

\begin{prop} \label {propisotropic2}
  Suppose that we have the standard situation where $C_2,\cdots, C_r,$ $ r\ge 2,$ are all isotropic and  $X={\rm Span}(C_2,\cdots,C_r)^\sharp=\bigcup_{i=2}^r C_i$. Assume that  $\dim X<n/2$. Then $C_1,\cdots,C_r$ are not the p-sets of a commutative Schur ring.
\end{prop}
 \noindent {\it Proof}  So let $p=\dim X <n/2$ and let   $a_1,\cdots ,a_p$ be a basis for $X$. Then we can extend this basis to a symplectic basis 
 $a_1,b_1,a_2,b_2,\cdots,a_{n/2}, b_{n/2}$ for $V$.  Here we have $b_i \in C_1,1\le i\le n/2$.
 Now we find that
 $$
 |\{x \in X:u \cdot b_1=0\}|=2^{p-1}-1,\qquad  |\{x \in X:u \cdot b_{n/2}=0\}|=2^{p}-1.
 $$
 It follows that $f_4(C_1)$ is not constant on $C_1$, a contradiction.\qed
 \medskip

 \begin{prop} \label {propisotropic3} Let $n\ge 6$. 
  Suppose that we have the standard situation where $C_2,\cdots, C_r$ are all isotropic, $r\ge 2$. Assume that
   $A={\rm Span}(C_2,\cdots,C_r)^\sharp=\bigcup_{i=2}^r C_i$ and  $\dim A=n/2$. Then $C_1,\cdots,C_r$ are not the p-sets of a commutative Schur ring.
\end{prop}
 \noindent {\it Proof}
So let   $e_1,\cdots ,e_{n/2}$ be a basis for $A$ where 
 $e_1,\cdots,e_{n}$  is a symplectic basis for $V$.
 Let $B={\rm Span}(\{e_{n/2+1},\cdots,e_{n}\})^\sharp$ so that
 $C_1=B \cup\left ( \cup_{b \in B} (b+A) \right )$.
 Define the multisets 
 $$
 X:=C_1^2 \cap (t_1t_2)^G,\quad Y:=C_1A\cap (t_1t_n)^G,\quad Z:=YA,
 $$
  so that $X, Y, Z \in \mathfrak S$.
 Now $Z$ has multiplicities $1,2,4$ and we let $Z_1$ be the set of elements of  $Z$ that have multiplicity $1$.
 We will show that $XZ_1 \ne Z_1X$. 
 
Note that $\a_1:=t_{1,2,3,n}t_{1,n-1}\in X,
  \a_2:=t_{2,n-1}t_{1,2,3}  \in Y,
  \a_3=\a_2t_{1,2} \in Z_1$.
 Let $\a=\a_1\cdot \a_3 \in X \cdot Z_1.$ 
 
 Let $V_6:={\rm Span}(e_1,e_2,e_3,e_{n-2},e_{n-1},e_n), W:={\rm Span} ( e_1+e_2+e_3+e_n, e_1+e_{n-1},
 e_2+e_{n-1},e_1+e_2+e_3,e_1+e_3  )={\rm Span}(e_1,e_2,e_3,e_{n-1},e_n)\le V_6$
 Then $\{u\a-u:u \in V_6\}=W$, showing that (i) $\a$ is not the product of less than five transvections, and (ii) if $t_at_bt_ut_vt_w=\a$, then $a,b,u,v,w $ is a basis for $W$. 
 
  Thus to consider whether $\a$ is in $XZ_1$ or $Z_1X$ we can  now restrict ourselves to considering products of transvections in
    $\mathcal T(W))$.
 This now reduces the problem  to checking the result in Sp$(W)\le G$. 
Specifically, 
one finds  that $\a$ has multiplicity $4$ in $XZ_1$, but multiplicity $0$ in $Z_1X$. This gives the result.
 \qed
 \medskip

 Thus we can now assume: $\Gamma^1(\mathcal T)$ is a star graph with $C_1$ the source, $f_1(C_2)>0, f_2(C_i)=0,i>1,$ and that $\bigcup_{i=3}^r C_i$ is isotropic.

                     \section {  Step 10:  The $f_1(C_2)>0$  case }

 \begin{prop}\label{propc1c2nn}
Consider the standard situation where $f_1(C_2)>0$. Let $X:=C_3\cup\cdots\cup C_r$.
Then $X$ is isotropic and  $X={\rm Span}(X)^\sharp$.
\end{prop}
 \noindent {\it Proof}   Each $C_i, i>2, $ is isotropic and $C_i \perp C_j$ for $2\le i \ne j\le r$. So $X$ is isotropic.
 
 Let $x \in {\rm Span}(X)^\sharp\setminus X$, so that $x \in C_1\cup C_2$
 also $x$ meets no elements of $C_2\cup\cdots \cup C_r$.
 If $x \in C_1$, then $x$ meet $2f_1(C_1)+f_2(C_1)$ elements of $C_1$ and meets no elements of $C_2\cup\cdots \cup C_r$.
 Thus $2f_1(C_1)+f_2(C_1)=2^{n-1}$.
 However there are   $y \in C_1$ that meet elements of  $C_2$ since $C_1\to C_2$ is an edge of $\Gamma^1(\mathcal T)$. Thus $y$  meets $2f_1(C_1)+f_2(C_1)=2^{n-1}$ elements of $C_1$ and meets at least one element of $C_2$. Thus $y$ meets more than $2^{n-1}$ elements of $\mathcal T$,  a contradiction. Thus $x \notin C_1$ and so  $x \in C_2$. 
Now $x \in  {\rm Span}(X)$ gives $x \perp C_2$. 
 Thus we cannot have  $x \in C_2$ since $f_1(C_2)>0$.
 So we have a contradiction.
\qed
\medskip

                     \begin{prop}\label{propfast1}
                     The standard situation where $C_1 \cap \left (  \bigcap_{j=2}^r C_j^\perp  \right)\ne \emptyset$ cannot happen.
                     \end{prop}
                      \noindent {\it Proof} Assume to the contrary  that   $b \in C_1 \cap \left (  \bigcap_{j=2}^r C_j^\perp  \right)$. Then the only elements of $\mathcal T$ that $b$ meets are  in $C_1$, and there are $2f_1(C_1)+f_2(C_1)$ of them; thus $2f_1(C_1)+f_2(C_1)=2^{n-1}$. 
                      However there is an edge of $\Gamma^1(C_1,\cdots,C_r)$ from $C_1$ to $C_2$. Thus there is some $c \in C_1$ that meets some element of $C_2$.
                      Since the $f_i$ are constant on $C_1$ we see that $c$ meets $2f_1(C_1)+f_2(C_1)$ elements of $C_1$; but it meets  at least one element of $C_2$.  Thus $b$ and $c$ meet a different number of elements of $\mathcal T$. This contradiction gives the result.\qed
                      \medskip
                      
                      Thus we may now assume 
                      \begin{align}\label{cbb1}
                      C_1 \cap \left (  \bigcap_{j=2}^r C_j^\perp  \right) =\emptyset.
                      \end{align}

Note that if none of the $C_i, 1\le i\le r$ is isotropic, then we must have $r=2$. So we now consider the situation where $C_3,\cdots, C_r$ are all isotropic and $r\ge 3$.

 Since $f_1(C_2)>0$ we may assume that $e_1,e_n \in C_2,$ so that $e_1+e_n \in C_2$. We may also assume that $e_2 \in C_3$.
 Let $a\in V, a \cdot e_2=1$. Then we must have $a \in C_1$, so that $e_{n-1}+e_2^\perp \subset C_1$. Let $W:={\rm Span}(e_1,e_n,e_2,e_{n-1})$. 
 Since $\Gamma^1(C_1,\cdots,C_r)$ is connected there 
 is $a \in C_1,a\cdot e_2=1$. Then $b:=a+e_2 \in C_1$ as $f_2(C_3)=0$.

 We recall that $\langle t_1,t_n,t_2\rangle  \le \langle C_2, C_3\rangle$ acts on each $C_i, 1\le i\le r$  and also acts on $\mathcal  T(W)$.
 Let $N_{12n}=N_G(\langle t_1,t_n,t_2\rangle).$
 One checks that there are two  orbits under the action of $N_{12n}$ on $\mathcal T$ for elements $a \in W$ where $a\cdot e_2=1$. Representatives for these orbits are $a\in \{e_{n-1},e_{n-1}+e_n\}$. Thus we can, without loss of generality, assume that 
 $a\in \{e_{n-1},e_{n-1}+e_n\}$. For such an $a$ we then take $b:=a+e_2$, so that $a,b \in C_1, a+b=e_2, a\cdot b=1$. 
 
 Then we have $a+e_2^\perp \subset C_1, e_1,e_n,e_1+e_n \in C_2, e_2 \in C_3$.  We also note that $a+e_2^\perp$ does not depend on the choice of $a$ if $a \cdot e_2=1$.
 This leaves three elements of $W^\sharp$ that are not known to be in $C_1\cup C_2 \cup C_3$, namely the elements of $R:=\{e_1+e_2,e_1+e_2+e_n,e_2+e_n\}$. 
 This leads us to do the following:
  Choose $R_1\subseteq R$ and $R_2 \subseteq R \setminus R_1$ and let
 $D_1=C_1\cup R_1, D_2=C_2 \cup R_2, D_3=C_3.$  We show

 \begin{lem}\label{lem123DDD} Assume the above and  $a \in \{e_{n-1},e_{n-1}+e_n\}.$
 Then the  only possibility for having $D_1,D_2,$ be the restriction to $W$ of p-sets $E_1,E_2$  of a commutative Schur ring is to have 
 $$
 D_1= a+e_2^\perp,\quad  D_2=\{e_1,e_n,e_1+e_n\} \cup R.
 $$
  Furthermore,  the set  $D_2$  is $\{u+v:u \in \{e_1,e_n,e_1+e_n\}, v \in \{0,e_2\}\}$.
 \end{lem}
 \noindent {\it Proof}
 This is similar to Proposition \ref {prop3cpts} and Lemma \ref{lemalpha}: let $X:=D_1D_2 \cap (t_1t_2)^G, Y=D_1^2 \cap (t_1t_n)^G, Z:=D_1Y$ and show that the only situation where $XZ$ is equal to $ZX$ is the situation described.\qed
 \medskip

 If $a \in \mathcal T$, then by $y_a$ we mean any element of $\mathcal T \setminus a^\perp$.

\begin{cor}\label{lemC2X} Assume that we have the standard situation where $f_1(C_2) >0$ and $r>2$. Then
 $x+C_2=C_2$ for all $x \in X=\bigcup_{i=3}^r C_i$. 
\end{cor}
\noindent {\it Proof}
Since $f_1(C_2)>0$, if $a \in C_2$, then there is $b \in C_2$ with $a\cdot b=1$ and $a,b,a+b \in C_2$.
 Now let $X:=C_3\cup \cdots \cup C_r$, an isotropic set orthogonal to $C_2$.
By  Lemma \ref {lem123DDD} in the above set-up  we see that for $x \in X$  we have $\{a+x,b+x,a+b+x\}\subseteq C_2$. 
 In particular $a+x \in C_2$ and so $x+C_2=C_2$ for all $x \in X$. \qed
 \medskip

 Thus 
 $
 \bigcup_{x \in X \cup\{0\}} \{x+e_1,x+e_n,x+e_1+e_n \} \subseteq C_2, \text { and }  \bigcup_{x \in X}  y_x+x^\perp \subseteq C_1
 $

 \begin{lem}\label{lemXC1}
 Let $X=C_3\cup \cdots \cup C_r$ (an isotropic set).  Let $x \in {\rm Span}( C_3\cup\cdots\cup C_r)$. Then $y_x+x^\perp \subseteq C_1$ and $
 {\rm Span}(X)^\sharp=X.$
 \end{lem}
 \noindent {\it Proof}
 Let $x \in {\rm Span}(X)$ and  
  $c\in y_x+x^\perp$, so that $c \cdot x=1$. Then $c \notin  X$ 
 as $X$ is isotropic, so that $c \in C_1\cup C_2$.
 If $c \in C_2$, then $c \cdot x \in \{c \cdot u:u \in X\} \ne \{0\}$; but $c \perp X$ gives a contradiction.
  Thus $c \in C_1$ and so $y_x+x^\perp \subseteq C_1$.
 
 Now if $x \in {\rm Span}( X) \setminus X$, then 
  $x \in C_1 \cup C_2$ and again $x \notin C_2$, so that
$x \in C_1$. 
   Now $f_2(C_1)>0$ means that there is some $c \in y_x+x^\perp\subseteq C_1, c\cdot x=1$,
    and $c+x\in C_2$. But $x \cdot (c+x)=1$, which   contradicts $X \perp C_2$.
 Thus ${\rm Span}(X)^\sharp = X$.
 \qed
  \medskip

\begin{prop} \label {propttrree}
  Suppose that we have the standard situation where $f_1(C_2)>0$ and  $C_3,\cdots, C_r$ are  isotropic, 
  $r\ge 3$. Let $X={\rm Span}  (C_3,\cdots,C_r).$ If $C_1,\cdots,C_r$ are   the p-sets of a commutative Schur ring, then
  either
  
  \noindent (i) $C_1 \cap C_2^\perp\cap X^\perp \ne \emptyset$; or
  
  \noindent (ii) $C_2^\perp=X.$
\end{prop}
 \noindent {\it Proof} So assume that (i) and (ii) are false,
  so that $C_1 \cap C_2^\perp\cap X^\perp =\emptyset$ and $C_2^\perp\ne X.$
But 
$X \subseteq C_2^\perp$ and   $C_2^\perp\ne X$ so we can choose $a \in C_2^\perp \setminus X$.
Now $a \in C_2^\perp$ and $f_1(C_2)>0$ means that $a \notin C_2$ and then $a \notin X$ gives 
$a \in C_1$.

Now we have $a \in C_1 \cap C_2^\perp$, so that 
$C_1 \cap C_2^\perp\cap X^\perp =\emptyset$ tells us that $a \notin X^\perp$.
 Thus there is some $b \in X$ such that $a \cdot b=1$.
 
 Let $c \in C_2$. Then $a \in C_2^\perp$ gives $a\cdot c=0$.
 By Corollary \ref {lemC2X} we have  $X+C_2=C_2$, so that $b+c \in C_2$. 
 However $a\cdot (b+c)=1$, contradicting $a \in C_2^\perp$. This contradiction completes the proof.
\qed
\medskip

By Lemma \ref {propfast1} we do not  have Lemma \ref{propttrree} (i), so we   now assume 
\begin{cor}\label{cor1iii}   Suppose that we have the standard situation where $f_1(C_2)>0$ and  $C_3,\cdots, C_r$ are  isotropic, $r\ge 3$. 
Then $C_2^\perp=X$.
\qed \end{cor}

\begin{prop} \label{proponenoniso}
Suppose that we have the standard situation where $f_1(C_2)>0, r\ge 3, n\ge 6,X={\rm Span}(X)^\sharp$ where $|X|>1$ and
$C_2^\perp=X$.
Then the $C_i$ cannot be the p-sets of a commutative Schur ring.
\end{prop}
\noindent{\it Proof} 
  Since  $X={\rm Span}(X)^\sharp$ where $|X|>1$ we can assume that $X={\rm Span}(a_1,a_2,\cdots,$ $a_k), 1<k<n/2,$ where 
 $a_1,b_1,\cdots,a_{n/2},b_{n/2}$ is a symplectic basis.
  Since $f_1(C_2)>0$ and  $f_2(C_2)=0$  we can assume: ${a_{k+1}},{b_{k+1}},{a_{k+1}+b_{k+1}} \in C_2$.
 Since $C_2+X=C_2$ we then see that 
  ${a_{k+1}+u},{b_{k+1}+u},{a_{k+1}+b_{k+1}+u} \in C_2$ for $u \in X$.
  We also know that $(y_{a_1}+a_1^\perp) \cup (y_{a_2}+a_2^\perp) \subseteq C_1$.  
 
 Let $W:={\rm Span}(a_1,b_1,a_2,b_2,a_{k+1},b_{k+1}), D_i=C_i \cap \mathcal T(W), i=1,2,$ and $D_X=\{a_1,a_2,a_1+a_2\}$. 
 Then $|D_1|=48$. Now  there are $12$ elements of $D_2$ of the form $u+v$ where $u \in \{a_{k+1},b_{k+1},a_{k+1}+b_{k+1} \}, v \in D_X\cup\{0\}$. Denote the set of these $12$ elements by $E_2$.
 Now let $E_1:=(b_1+a_1^\perp) \cup (b_2+a_2^\perp).$ Then $E_1,E_2,D_X$ are disjoint and we have $E_1\subset D_1, E_2\subset D_2$. So we obtain
   $|E_1|+|E_2|+|E_X|\ge 48+12+3=63=|W^\sharp|$, showing that $D_1=E_1,D_2=E_2,D_X$ is a partition of $W^\sharp$.  

Now let $A:=(C_1^3)_1,$ the set of elements of $C_1^3$ that have multiplicity $1$ in $C_1^3$, and let
$B:=(C_1C_2C_1)_1$  the set of elements of $C_1C_2C_1$ that have multiplicity $1$ in $C_1C_2C_1$. Then $A,B \in \mathfrak S$.
 Now let 
 \begin{align*}&
 \a_1:=t_{b_2+b_1}t_{a_1+a_{k+1}+b_2+b_1}t_{a_1+b_{k+1}+b_2+b_1},\\
& \a_2=t_{a_2+a_{k+1}+b_{k+1}+b_2}t_{a_2+a_{k+1}+b_{k+1}},t_{a_1+a_2+b_2},
\end{align*}
 and note that $\a_1 \in A, \a_2 \in B$. 
 Let $\a=\a_1\a_2\in AB$. We note that $\{v\a-v:v \in W\}=W$, showing that the minimal number of transvections giving $\a$ as a product is six, and that any such transvections must be in $\mathcal T(W)$, and must span $W$. Thus to obtain $\a$ as a product $\b_1\b_2, \b_1 \in A,\b_2 \in B$ (or where $\b_1 \in B, \b_2 \in A)$ we only need to consider products of transvections in $\mathcal T(W)$.
 
Now the multiplicity of $\a$ in the multi-set $AB$ is $6$, while 
  the multiplicity of $\a$ in  $BA$ is $3$, so that $AB \ne BA$.\qed
\medskip

%%CHECKED

\begin{cor}\label
{cor32111}
Suppose that we have the standard situation where $f_1(C_2)>0, r\ge 3, n\ge 6,X={\rm Span}(X)^\sharp$. Then $r=3$ and we can assume that $|X|=|C_3|=1$  and 
$C_2^\perp=X$. \qed
  \end{cor}

                    \section {    Step 11 }
 
 \subsection {The $r=2$ case} Assume that we have the standard situation with $r=2.$
                     \begin{prop} \label{propconn} 
              If $ f_1(C_2)>0, f_2(C_2)=0, f_2(C_1)>0$, then 
                     
\noindent                     (i) $\Gamma^1(C_2)$ is connected, and 
                     
               \noindent      (ii) $\langle C_2\rangle$ is a 3-transposition subgroup of ${\rm Sp}(n,2)$ (where $C_2$ is the set of involutions 
               of  $\langle C_2\rangle$ 
               determing the $3$-transposition group structure)  that acts transitively on $C_2$ by conjugation and also acts on $C_1$ by conjugation.
                     \end{prop}
                                \noindent {\it Proof} (i) Suppose that $\Gamma^1(C_2)$ is not connected and let $\Gamma_{2i}, 1\le i\le s,$ be the components of  $\Gamma^1(C_2)$. 
                                Then  $f_{1}(\Gamma_{2i})=f_1(C_2)>0, 1\le i\le s$.
                                  \smallskip
                               
                  \noindent             {\bf CASE 1:} $s\ge 3$. Here we choose $b_i \in \Gamma_{2i}$. Then $b_i$ does not meet any elements of $ \Gamma_{2j}, j \ne i$. Then  $b_1+b_2, b_1+b_2+b_3 \in C_1$ (as these element meet 
                                $ \Gamma_{21}$ and $ \Gamma_{22}$). But they meet the same number of elements of $C_1, \Gamma_{21}$ and $  \Gamma_{22}$, but not the same number of elements of $ \Gamma_{23}$ (and they meet no elements of $ \Gamma_{2i}, i>3$). This contradiction gives the result in this case. So we are left to consider:
                                \smallskip
                                
                 \noindent              {\bf CASE 2:}   $s=2$. 
                                We first claim that ${\rm Span}( \Gamma_{2i})^\sharp= \Gamma_{2i}, i=1,2$. 
                                To this end assume that ${\rm Span}( \Gamma_{21})^\sharp  \ne  \Gamma_{21}$ and let $b_0 \in {\rm Span}( \Gamma_{21})^\sharp  \setminus  \Gamma_{21}$. 
                                Now choose $b_2 \in  \Gamma_{22}.$  Then $b_0, b_0+b_2 \in C_1$ and   they both meet the same number of elements of $C_1$ and of $ \Gamma_{21}$, however $b_0$ meets zero elements of $ \Gamma_{22}$, while $b_0+b_2$ meets $2f_1(C_2)>0$ elements of $ \Gamma_{22}$. This contradiction gives the claim. It then follows that each $\Gamma_{2i}$ is a (non-degenerate) symplectic subspace of $V$, since $0<f_1(C_2)=f_1(\Gamma_{2i})$ and so there are no  isotropic elements of  $\Gamma_{2i}, i=1,2$.
                                
            So we now have $ \Gamma_{21}=V_u^\sharp,  \Gamma_{22}=V_v^\sharp, u+v\le n$, where $V_k$ is a symplectic space of dimension $k$ and $V_u\perp V_v$. Also $V_u \cap V_v=\{0\}$.
            
            We next claim that ${\rm Span} (C_2)=V$. If not, then from the above there is a symplectic subspace $V_k, k>1,$ such that
            $V=V_u\oplus V_v\oplus V_k$. Let $b_0 \in V_k^\sharp, b_i \in \Gamma_{2i},i=1,2$.  
            Then $b_0+b_1,b_0+b_1+b_2\in C_1$ again meet different numbers of elements of $V$. Thus $V=V_u \oplus V_v.$ [This is the end of the proof of Proposition \ref {propconn3} - see below.]

            Since $n\ge 6$ and $u,v >0$ we can assume that $u\ge 4, v \ge 2$. 
            Let 
            $$W:={\rm Span}(e_1,e_2,e_3,e_{n-2},e_{n-1},e_n), U:={\rm Span}(e_1+e_6,e_2,e_3+e_6,e_4+e_6,e_5+e_6)
            $$
             and 
       assume that ${\rm Span} (e_1,e_2,e_{n-1},e_n) \subseteq V_u, {\rm Span} (e_3,e_{n-2}) \subseteq V_v$. Then 
            $W \setminus ({\rm Span} (e_1,e_2,e_{n-1},e_n) \cup  {\rm Span} (e_3,e_{n-2}))\subset C_1$.
            Define the following multisets:  
            $$
            X:=(C_1\cdot C_2) \cap (t_1t_2)^G,\quad  Y:=C_1^2\cap (t_1t_n)^G, \quad Z:=C_1 \cdot Y.
            $$ Then 
            $X,Y,Z \in \mathfrak S$. 
            Let $Z_i$ be the set of those elements of the set $Z$ that have multiplicity $i\in \mathbb Z$ in the multiset $Z$.
           We will show that $X   Z_1 \ne Z_1    X$. 
            Let 
            $$
          z:=t_{n-2,n}t_{2,n-2,n-1}t_{n-2,n-1} \in Z_1, \quad x:=t_{2,3,n}t_{1,2,n-1} \in X,\quad \a:=zx\in Z_1X.
            $$
            Then  $\a$ is a product of $5$ transvections each  in $\mathcal T(U)$, and  that if $\a=t_ut_vt_wt_xt_y$, then ${\rm Span}(u,v,w,x,y)=U$. Thus to find such $u,v,w,x,y$ we need only look in $U$.
            
            Now  $\mathcal T(U)$ is partitioned as  $(C_1 \cap \mathcal T(U)) \cup (C_2 \cap \mathcal T(U))$ and it is  easy to check that $\a$ occurs in $Z_1X$ with multiplicity  $5$ while $\a$ occurs with multiplicity  $3$ in $XZ_1$. Thus $XZ_1 \ne Z_1X$ and this concludes the proof of Proposition \ref {propconn} (i).
            
            By (i)  $\Gamma^1(C_2)$ is connected, and  if $t_a,t_b \in C_2, a\cdot b=1$, then $t_a^{t_bt_a}=t_b.$ Thus   elements of $ C_2$ are  conjugate by elements of $\langle C_2\rangle$. 
            Now   $C_2$ is fixed by the conjugation action of any element of $\langle C_2\rangle$. Thus $C_2$ is a conjugacy class in $\langle C_2\rangle.$
            This  shows that $\langle C_2\rangle$ is a 3-transposition subgroup of ${\rm Sp}(n,2)$ and the rest follows.
            \qed
            \medskip

               \begin{lem}\label{lemconnC2} Assume:  $r=2, f_2(C_2)=0, f_1(C_2)>0, f_2(C_1)>0,$ and   $\Gamma^1(C_2)$ is connected. Then 
               ${\rm Span}(C_2)=V$ and $\langle C_2\rangle$ is an irreducible $3$-transposition group.
               \end{lem}
                          \noindent {\it Proof} 
                          So assume that $0\ne x \in {\rm Span}(C_2)^\perp$. Since $f_1(C_2)>0$ we have  $x \notin C_2$, so that $x \in C_1$. Thus $x \in C_1 \cap (C_2^\perp)=\{0\}$,  contradicting $f_2(C_1)>0$.
                          Thus ${\rm Span}(C_2)$ is a symplectic subspace and if ${\rm Span}(C_2) \ne V$, then there is $0\ne x\in C_2^\perp$ and again we have $x \in C_1 \cap (C_2^\perp)$.\qed
                          \medskip

                           \subsection {The $r=3$ case}   Assume  the standard situation where $\mathcal T=C_1\cup C_2\cup C_3, C_3=\{e_1\},
                           f_1(C_2)>0, f_2(C_2)=0, f_2(C_1)>0$.
                          
                     \begin{prop} \label{propconn3} 
             Then 
                     
\noindent                     (i) $\Gamma^1(C_2)$ is connected, and 
                     
               \noindent      (ii) $\langle C_2\rangle$ is a 3-transposition subgroup of $Sp(n,2)$ (where $C_2$ is the set of involutions)  that acts transitively on $C_2$ by conjugation and also acts on $C_1$ by conjugation.
                     \end{prop}
                                \noindent {\it Proof} (i)
  The proof   follows the proof of   Proposition \ref {propconn} down to where it is shown that
  $V_u\oplus V_v=V$. This   gives a contradiction as then $C_3\perp C_2$ and ${\rm Span}(C_2)=V$.  (ii) follows as in  Proposition \ref {propconn}.          \qed

 %XXXX                             

\section{Step 12:   Relations}

\subsection {Relations  when  $r=2$} So  assume that $r=2, n\ge 6, \mathcal T=C_1 \cup C_2$ where $f_2(C_2)=0, f_2(C_1)>0$.
Let $a_i=f_i(C_1), b_i=f_i(C_2), 1\le i\le 4$. 
The goal will be to find polynomial relations among $a_1,\cdots,a_4,b_1,\cdots,b_4$ and some other indeterminates.
We certainly have (see Eq. (\ref{eq123456}))
\begin{align}
&
\label{ir1} a_3=2a_1+a_2+a_4+1;\\
&
\label{ir2} b_3=2b_1+b_2+b_4+1;\\
&
\label{ir3} a_3+b_3=2^n-1;\\
&b_2=0.
\end{align}

     For $0 \ne a \in V, e \in \mathbb F_2,$ and $C \subseteq \mathcal T$ we define
     $B_a^e(C)=\sum_{t_b \in C,a \cdot b=e} t_b.$
     
     \begin{lem} \label{lemabc} 
     Let $C_1,C_2 \subseteq \mathcal T$ be disjoint  $\mathfrak S$-sets.
Then 

\noindent (i) $\sum_{t_a \in C_1} B_a^0(C_2) \in \mathfrak S$. 

\noindent (ii) If $C_2$ is a p-set, then  $\sum_{t_a \in C_1} B_a^0(C_2) =\lambda C_2$ for some integer $\lambda\ge 0$. 

 \noindent (iii) If  $\mathcal T= C_1 \cup C_2$ and  $t_a \in C_1$, then $|B_a^0(C_2)|=2^{n-1}-2-f_4(C_1).$

\end{lem}
\noindent {\it Proof}
Now the $(t_1t_2)^G$ part of $\overline {C_1}\,\overline {C_2}$ is $U:=\sum_{t_a \in C_1} t_aB_a^0(C_2)$ since $C_1 \cap C_2=\emptyset$. 
So
$$
\overline{C_1} \cdot U=\sum_{t_{a'},t_{a} \in C_1} t_{a'}t_a B_a^0(C_2).
$$

Let $t_{a'}t_a t_b\in \overline{C_1} \cdot U$,   where $t_{a'},t_a\in C_1, t_b \in B_a^0(C_2)$, so that $a \cdot b=0$. Since $a,a' \ne b$ Lemma \ref {lem43} shows that the only way that $t_{a'}t_a t_b \in \mathcal T$ is if $a=a'$.

Thus we see that the $\mathcal T$ part of this sum is $\sum_{T_{a} \in C_1}  B_a^0(C_2)$, which is thus in $\mathfrak S$. This gives (i) and (ii) follows. 

(iii) Note that $|a^\perp |=2^{n-1}$. Since there are $f_4(C_1)$ elements $t_b \in C_1, b \ne a,$ with $b \cdot a=0$, there are (since $\mathcal T=C_1 \cup C_2$) $2^{n-1}-2-f_4$ elements $t_c \in C_2$ with $c \cdot a=0.$  This gives $|B_a^0(C_2)|=2^{n-1}-2-f_4(C_1)$.\qed\medskip

\begin{cor} \label {corlemabc}
If $\mathcal T= C_1 \cup C_2$ where $C_1,C_2$, then  there are $\l_1,\l_2\in \mathbb Z^{\ge 0}$ with 
\begin{align} &
 \label{ir5}  a_3(2^{n-1}-2-a_4)=\l_1 b_3;\\
 &
  \label{ir6}  b_3(2^{n-1}-2-b_4)=\l_2 a_3.\qed
\end{align}
\end{cor}

%We introduce the following notation in the situation where  $\mathcal T=C_1\cup C_2$ is the union of two p-sets. Let $a_i=f_i(C_1),b_i=f_i(C_2),i=1,2,3,4.$

\begin{lem}\label{lemsizes2}
As a function of the $a_{i},b_j$ the number  of zero-triangles   in $\mathcal T$ is
  \begin{align}\label{ir7}
    \frac 1 3 a_1a_3+\frac 1 2 a_2a_3+\frac 1 3 b_1b_3+\frac 1 2 b_2b_3=\frac 1 3 (2^n-1)2^{n-2}.
    \end{align}
\end{lem} 
\noindent {\it Proof} Any $t_a \in C_1$  is in $f_1(C_1)=a_1$ zero-triangles contained in $C_1$. Summing over all $a_3$ elements of $C_1$ we get a contribution of $a_{1}a_3/3 $ to the number of zero-triangles in $\mathcal T$.
Further $t_a\in C_1$  is in $a_2$ zero triangles 
$\{t_a,t_b,t_{a+b}\}$ for some $t_b \in C_1$ and  $t_{a+b} \in C_2$.
 Summing over all $a_3$ elements of $C_1$ we obtain a contribution of $a_{2}a_3/2 $ to the number of zero-triangles in $\mathcal T$.
Now we do the same argument for $t_a \in C_2$, noting that counting in this way we never count any zero-triangle twice. Thus the total number is 
 $
a_{1}a_3/3 + a_{2}a_3/2+
b_{1}b_3/3 + b_{2}b_3/2.
$ We now use Lemma \ref {lemsizes}. \qed\medskip

 Assume $\mathcal T=C_1 \cup C_2 \cup\cdots\cup C_r$. 
 Let $\Sigma(C_2,C_1,a), a \in C_2,$ be the graph with vertices $C_1$ and  edges $\{v_1,v_2\}\subset C_1$ if $v_1\cdot v_2=1$ and $v_1+v_2=a$.
 For $S \subseteq C_2$ let
 $$
 \Sigma(C_2,C_1,S)= \bigcup_{a \in S} \Sigma(C_2,C_1,a)\text { and }  
  \Sigma(C_2,C_1)=   \Sigma(C_2,C_1,C_2)
 $$
 
 We first note:  if we have an edge in $\Sigma(C_2,C_1,a)$, then $f_2(C_1) \ne 0$.
  Suppose   $\mathcal T=C_1 \cup C_2$ and  $f_2(C_1)\ne 0$. Then for   $b \in C_1$ there are $f_2(C_1) $ elements $a+b \in C_1$ with $a\cdot b=1$ and $a \in C_2$. This gives $f_2(C_1)$ elements $a\in C_2$ for each $b \in C_1$. 
  
  Conversely, for $a \in C_2$, if there is $b 
  \in C_1$ with $a\cdot b=1, a+b \in C_1$, then this $b\in C_1$ is accounted for as in the above paragraph.

Fix $a \in C_2$.  Now the number 
of $b \in \mathcal T$ with $a\cdot b=1$ is $2^{n-1}$ of which $2f_1+f_2$ are in $C_2$. Thus there are $2^{n-1}-2f_1-f_2$ such elements $b$   in $C_1$.  
 There are $a_2$ of the $b \in C_1$ with $a+b \in C_2$. This leaves $2^{n-1}-2a_1-2a_2$  elements $b \in C_2$ where $a+b \in C_1$.  These $b$s come in pairs $b,a+b$, each such pair giving an edge of $\Sigma(C_2,C_1,a)$.  Thus the number of edges of $\Sigma(C_2,C_1,a)$ is $b_3(2^{n-1}-2b_1-2b_2)/2$.
 But the number of such edges is also $f_3(C_1)f_2(C_1)/2=a_2a_3/2$. Thus we have
  \begin{align}\label{ir8}
 &a_2a_3=b_3(2^{n-1}-2b_1-2b_2) 
\end{align}
                             
\subsection {Relations for $r=3$}
So assume that $r=3, n\ge 6$ and that $\mathcal T=C_1 \cup C_2\cup C_3, C_3=\{e_1\}$ where $f_2(C_2)=0, f_2(C_1)>0$.
We also have $e_n+e_1^\perp \subseteq C_1, C_2+e_1=C_2.$
Let $a_i=f_i(C_1), b_i=f_i(C_2), 1\le i\le 4$. 
We again find polynomial relations among $a_1,\cdots,a_4,b_1,\cdots,b_4,\l_1,\l_2$.
We certainly have
\begin{align}
&
\label{ir13} a_3=2a_1+a_2+a_4+1;\\
&
\label{ir23} b_3=2b_1+b_2+b_4+1;\\
&
\label{ir33} a_3+b_3=2^n-2;\\
&b_2=0.
\end{align}

An argument similar to Corollary \ref {corlemabc} gives:
\begin{lem} \label {corlemabc3}
If $\mathcal T= C_1 \cup C_2$ where $C_1,C_2$ are p-sets, then  
$
f_3(C_1)(2^{n-1}-2-f_4(C_1))=\l f_3(C_2)
$
for some integral $\l\ge 0$. Thus
\begin{align} &
 \label{ir53}  a_3(2^{n-1}-2-a_4)=\l_1 b_3;\\
 &
  \label{ir63}  b_3(2^{n-1}-3-b_4)=\l_2 a_3.\qed
\end{align}
\end{lem}

\section {Symmetric group cases}

We now consider the situations where $\langle C_2\rangle$ is a symmetric group $S_{n+1}$ or $S_{n+2}$.
             
Note that for $r=2,3$:   if we know $a_3$ or $b_3$ as a function of $n$, then we know $a_3$ and $b_3$ as a function of $n$ (since $a_3+b_3=2^n-1$ or $a_3+b_3=2^n-2$). Then the quadratic equations  Eqs. (\ref {ir1})-(\ref {ir8}) in the   $a_i,b_i, 1\le i\le 4,\l_1,\l_2$ 
(or Eqs. (\ref {ir13})-(\ref {ir63}))
become linear equations with coefficients that are functions of $n$.

\subsection {Symmetric group $S_{n+1}$ for $r=2$}
So assume that $\mathcal T=C_1\cup C_2$ and 
 note that if $\langle C_2\rangle=S_{n+1}$, then 
$
b_3=\binom {n+1} 2 \text { and  } b_4=\binom {n-1}2.
$
Using Eq (\ref {ir3})  gives 
$
a_3=2^n-1-\frac 1 2 n^2-\frac 1 2 n.
$
Using $b_2=0$ and  these values for $a_3,b_3, b_4$ Eqs. (\ref {ir1}) -(\ref {ir8}) become
\begin{align*}
&2a_1 + a_2 + a_4 - 2^n + \frac 1 2 n^2 + \frac 1 2 n + 2=0;\quad 
    2b_1  - 2n + 2=0,\\&
  \left   (\frac 1 3 2^n - \frac 1 6 n^2 - \frac 1 6 n - \frac 1 3  \right ) a_1 +   \left  ( 2^{n-1} - \frac 1 4 n^2 - \frac 1 4 n - \frac 1 2\right ) a_2 + 
  \left          (\frac 1 6 n^2 + \frac 1 6 n\right ) b_1   - \frac 1 {12} 2^{2n} + \frac 1 {12} 2^n=0,\\&
     \left   (-2^n + \frac 1 2 n^2 + \frac 1 2 n + 1\right ) a_4 +  \left   (-\frac 1 2 n^2 - \frac 1 2 n \right ) \l_1 + \frac 1 2 2^{2n} - \frac 1 4 2^n n^2
       - \frac 1 4 2^n n - \frac 5 2  2^n + n^2 + n + 2=0,\\&      \left  (-2^n + \frac 1 2 n^2 + \frac 1 2 n + 1\right ) \l_2 + \frac 1 4 2^n n^2 + \frac 1 4 2^n n - \frac 1 4 n^4 + \frac 1 2 n^3 - 
        \frac 3 4  n^2 - \frac 3 2  n=0,\\&
    \left    (- 2^{n+1} + n^2 + n + 2\right ) a_1 +   \left  (-2^{n+1} + n^2 + n + 2\right ) a_2  + 
        2^{2n-1} - \frac 1 4 2^n n^2 - \frac 1 4 2^n n - \frac 1 2 2^n=0,\\&
      \left  (2^n - \frac 1 2 n^2 - \frac 1 2 n - 1 \right ) a_2 +   \left  (-\frac 1 2 n^2 - \frac 1 2 n \right ) \mu=0.
    \end{align*}
We solve these linear equations in the variables $a_1,a_2,a_4,b_1,\l_1,\l_2,\mu$ to get 
 $a_1=N/D$ where
$$
N=2^{2n-2} - 3\cdot 2^{n-3}n^2 - 3\cdot 2^{n-3}n - 2^{n-2} + n^3 - n;\quad
D=2^n - \frac 1 2 n^2 - \frac 1 2 n - 1.
$$
One finds that $N/D \notin  \mathbb Z$ if $n\ge 8$, a contradiction.\qed

 \subsection {Symmetric group $S_{n+2}$ for $r=2$}
 The method is  the same as  for the $S_{n+1}$   case above, except that $ 
 b_3=\binom {n+2} 2, b_4=\binom {n}2$.
 One finds that $a_1=N/D$ where
 $$
 N=2^{2n-2} - 3\cdot 2^{n-3} n^2 - 9\cdot 2^{n-3}n - 2^n + n^3 + 3\cdot n^2 + 2n,\quad D=2^n - \frac 1 2 n^2 - \frac 3 2 n - 2.
 $$
 One then  checks that $N/D$ is not an integer for all $n\ge 14.$

\section {Orthogonal Groups  for $r=2$}

We recall the definitions of quadratic forms, their Arf invariants and  their Orthogonal groups.
A {\it quadratic form on $V$} (a vector space over $\mathbb F_2$) is a function $Q:V  \to \mathbb F_2$ satisfying
$Q(ax+by)=aQ(x)+bQ(y) +(ab)(x \cdot y)$ for $a,b \in \mathbb F_2, x,y \in V$. We note that then $x \cdot y=Q(x+y)+Q(x)+Q(y)$.
Put $n=2m$. Then  \cite [pp. 197-199]   {dic}
 the quadratic form $Q$ is equivalent to $Q_0$ or $Q_1$ where
$$Q_\a\left (\sum_{i=1}^n x_ie_i\right )=( \alpha x_1+x_1x_n+\alpha x_n)+x_2x_{n-1}+x_3x_{n-2}+\cdots +x_mx_{m+1}.
$$
 Here $\a \in \mathbb F_2$. Further  $Q_0$ and $Q_1$  are not equivalent and these give the two types of orthogonal group.

The {\it Arf invariant} of the quadratic form (which takes values in $\mathbb F_2$)  has a number of definitions: it is $0$ if the quadratic form is equivalent to a direct sum of copies of the binary form 
$xy$, and it is $1$ if the form is a direct sum of 
$
x^{2}+xy+y^{2}$ with a number of copies of 
$
xy.$

The Arf invariant  is also  the value which is assumed most often by the quadratic form \cite {Br}.
 Another characterization: $Q$  has Arf invariant $0$ if and only if the underlying $n$-dimensional vector space over the field $\mathbb F_2$ has a $n/2$-dimensional subspace on which $Q$ is identically $0$.
 
If $a_1,b_1,\cdots,a_{n/2},b_{n/2}$ is a symplectic basis then also: ${\rm Arf}(Q)=\sum_{i=1}^{n/2} Q(a_i)Q(b_i)$; here it is known that this is independent of the choice of symplectic basis.

Thus there are two types of quadratic form, as determined by their Arf invariant. These give rise to two Orthogonal groups
${\rm SO}^{\pm}(n,2)$. The form $Q_0$ corresponds to the Orthogonal group ${\rm SO}^-(2m,1)$ 
while
the form $Q_1$ corresponds to  ${\rm SO}^+(2m,1)$.

Given a basis $B=\{b_1,\cdots,b_n\}$ for $V$ a quadratic form $Q_B$ can be defined  using $B$: just take $Q(b_i)=1$ for all $i\le n$ and extend using the above equation.
If $H=\langle B\rangle$ and the graph $\Gamma^1(B)$ is connected, then $H$ acts transitively on $B$ and so on the orbit $(B)H$, this being the conjugacy class of symplectic transvections in $H$. 
Further, the  quadratic form $Q_B$ is invariant under the action of $H$: $Q((a)h)=Q(a)$ for all $a \in V, h \in H$, and $H$ is the corresponding orthogonal group.

\begin{prop}\label{propOT}
The group  ${\rm SO}^+(2m,2)$ has $2^{m-1}(2^m-1)$ symplectic transvections, while 
${\rm SO}^-(2m,2)$  has $2^{m-1}(2^m+1)$ symplectic transvections.
\end{prop}
\noindent {\it Proof}
This will be by simultaneous  induction on $m\ge 2$, the case $m=2$  being easy to check.
Thus we assume that 
 the group ${\rm SO}^+(2m,2)$ has $2^{m-1}(2^m-1)$ symplectic transvections, while the group 
${\rm SO}^-(2m,2)$  has $2^{m-1}(2^m+1)$ symplectic transvections. We then consider the case  $m+1$.

We write the form $\sum_{i=1}^{m+1} x_iy_i$ as $\left (\sum_{i=1}^{m} x_iy_i\right )+(x_{m+1}y_{m+1})$.
We note that $x_{m+1}y_{m+1}$ is zero for exactly three choices of  $x_{m+1},y_{m+1} \in \mathbb F_2$ (and is $1$ otherwise).
 
 By induction, the number of $(x_1,y_1,\cdots,x_m,y_m)$ with $\sum_i x_iy_i=0$ is $2^{m-1}(2^m-1)$ while the number with 
  $\sum_i x_iy_i=1$ is $2^{m-1}(2^m-1)$, so the number of $(x_1,y_1,\cdots,\break x_{m+1},y_{m+1})$ with $\sum_i x_iy_i=0$ 
  is
 $$
  3\cdot 2^{m-1}(2^m-1)+1\cdot 2^{m-1}(2^m+1)=2^{m+1-1}(2^{m+1}-1),
 $$ 
 as required. 
 For the other case the sum that we get is
 $
 3\cdot 2^{m-1}(2^m+1)+1\cdot 2^{m-1}(2^m-1)=
2^{m+1-1}(2^{m+1}+1),
 $ 
 as required. \qed

\begin{prop}\label{propOTf4}  
For $H={\rm SO}^\pm (n,2),n\ge 8, $ we have  $f_4(C_2)=2^{n-2}-1$.
\end{prop}
\noindent {\it Proof}
We  induct on $n\ge 8$
We   use a set $B=B_n$ of transvection generators for $H=H_n\in \{{\rm SO}^+(n,2),{\rm SO}^-(n,2) \}$ given by \cite {SS}. See  Figure 1, where there are two graphs $J_{1}, J_{2}$ shown. 
We let $J_{n,i}$ be the subgraph of $J_i, i=1,2,$ on vertices labeled $1,2,\cdots,n$.
The relationship between $B$ and the graphs is that $\Gamma^1(B)$ is one of the graphs $J_{n,1}, J_{n,2}$; in each case $B$ is a basis for $V$.
 For example, the graph for ${\rm SO}^+(n,2)$ is the graph $J_{n,i}$,  where vertex $n$  is labelled 
${\rm O}^+$ in $J_i$. 
Vertex $j$ of the  graph $J_{n,i}$ corresponds to $b_j\in B$ if $B=\{b_1,\cdots,b_n\}$. Each graph $J_i,i=1,2,$ is  extended to the right with vertex labelings of  period $8$. 
  Let $Q_B$ be the corresponding quadratic form.

\begin{figure}
     \centering
   {
          \includegraphics[width=.85\textwidth]{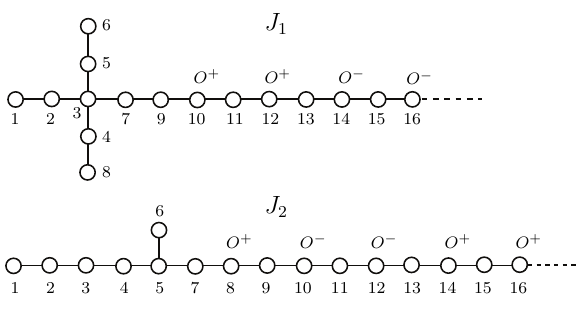}}
\caption{}
\end{figure}

For $H \le G$ and  any $t \in \mathcal T(H)$ let $B(t)=(B(t)_i)_{i=1}^n$ be the coordinates of $t$ relative to the basis $B=\{b_1,\cdots,b_n\}$.
Now for $\varepsilon , \delta \in \mathbb F_2$ we define:
$$
n_{\varepsilon,\delta}= |\{t\in \mathcal T(H): Q_B(t)=\varepsilon, B(t)_n=\delta, b_1 \cdot t=0\}|.
$$
Note that the condition $b_1\cdot t=0$ is equivalent to $B(t)_2=0$. We also note that for each $v \in V$ the value $Q_B(v)$ is just the number of components (mod $2$) of the subgraph of $\Gamma^1(B)$ determined by $v$ i.e. determined by the non-zero coordinates of $B(v)$.
Also  $t\in \mathcal T(H)$ determines a pair $\pi(t):=[\varepsilon,\delta]$ where $Q_B(t)=\varepsilon, B(t)_n=\delta.$

One checks the $n=8$ case to start the induction \cite {ma}. Now to get from the $n$ case to the $n+2$ case we just add two vertices ($n+1$ and $n+2$) to the right of vertex $n$. 
We want to find linear equations relating the $(n+2)_{\varepsilon,\delta}$ to the $n_{\varepsilon,'\delta'}$.
Now any $u \in V_{n+2}$ can be written as $u=v+w$ where $v \in V_n$ and 
$w=\l_1 b_{n+1}+\l_2b_{n+2}$. So, for fixed $v,$ there are $4$ cases. 
For example, if $t_v\in \mathcal T(H_n)$ has $\pi(t_v)=[0,1]$ and we add $0b_{n+1}+1b_{n+2}$, then $\pi(t_{v+b_{n+2}})=[1,1]$.
Considering all such possibilities gives the following system of equations:
\begin{align*}
& (n+2)_{0,0}=n_{0,0}+n_{1,0}+2n_{0,1}+1;\\
& (n+2)_{1,0}=n_{0,0}+n_{1,0}+2n_{1,1}+1;\\
&(n+2)_{0,1}=2n_{1,0}+n_{0,1}+n_{1,1}+2;\\
&(n+2)_{1,1}=2n_{0,0}+n_{0,1}+n_{1,1}+2.
\end{align*}
As a matrix equation this is:
\begin{align*}&
\begin{bmatrix}
1&1&2&0&1\\
1&1&0&2&1\\
0&2&1&1&2\\
2&0&1&1&2\\
0&0&0&0&1
\end{bmatrix}
\begin{bmatrix}
n_{0,0}\\
n_{1,0}\\
n_{0,1}\\
n_{1,1}\\
1
\end{bmatrix} =\begin{bmatrix}
(n+2)_{0,0}\\
(n+2)_{1,0}\\
(n+2)_{0,1}\\
(n+2)_{1,1}\\
1
\end{bmatrix} 
\end{align*} 
The extra $1,1,2,2$ (in the last column of the $5 \times 5$ matrix) that are added come from considering the cases $t_v=t_0=id$ and $t_v=t_{b_1}$.
Denote the above $5 \times 5$ matrix by $K$. Then $K$ has eigen-values $0,4,1,2i,-2i,i=\sqrt{-1}$, and so is diagonalizable with diagonalizing matrix  \cite {ma}
$$Q:=
\begin{bmatrix}
     4    &   4 &       -4  &      -4  &       8\\
     60    &    60    &    60   &     60   &    120\\
      0 &        0 &        0 &        0 &       15\\
 20i&   -20i    &    20   &    -20   &      0\\
-20i &   20i  &      20 &      -20 &        0
\end{bmatrix}
$$ Here $QKQ^{-1}={\rm diag}[0,4,1,2i,-2i]$. 
One then shows that for $m\ge 0$ we have: 
$$
K^m\begin{bmatrix}
31,31,32,32,
1
\end{bmatrix} ^T
=\begin{bmatrix}
2^{2m+5}-1,2^{2m+5}-1,2^{2m+5},2^{2m+5},
1
\end{bmatrix}^T.
$$
With this,  Proposition \ref {propOTf4}  follows as the number required is either $n_{0,0}+n_{0,1}$ or $n_{1,0}+n_{1,1}$, which in either case is $2^{n-2}-1$.\.
\qed

 \subsection {Orthogonal group ${\rm SO}^+(n,2)$ for $r=2$}

Let $H={\rm SO}^+(n,2), n=2m$.
We now have  $b_4=f_4(C_2)=2^{n-2}-1$  (Proposition \ref {propOTf4})  and $b_3=2^{m-1}(2^{m}-1)$ (Proposition \ref {propOT}). We also have  Eqs (\ref {ir1}), (\ref {ir2}) and (\ref {ir3}). These five equations are  linear and we solve them to get
values for $a_1,a_3,b_1,b_2(=0),b_3,b_4,\l_1,$ $\l_2$. Substituting these values into the remaining equations we obtain:
\begin{thm}\label{thmOpos} Let $H={\rm SO}^+(n,2), n=2m$.
Then  $a_1,\cdots,b_4,\l_1,\l_2$ are
\begin{align*}
& a_1=0,&&
    a_2 = 2^{n-2},\\
&    a_3 = 2^{n-1} + 2^{m-1} - 1,&&
    a_4 = 2^{n-2} + 2^{m-1}  - 2,\\
 &   b_1 = 2^{n-3} - 2^{m-2} ,&&
    b_2=0,\\
    &b_3 = 2^{n-1} - 2^{m-1} ,&&
    b_4 = 2^{n-2} -1,\\
  &  \l_1 = 2^{n-2}- 1,&&
    \l_2 =2^{n-2} - 2^{m-1} .\qed 
\end{align*}
\end{thm}

 \subsection {Orthogonal group ${\rm SO}^-(n,2)$ for $r=2$}
 A similar calculation, again using Propositions \ref {propOTf4} and \ref {propOT}, deals with the case $H={\rm SO}^-(n,2)$:

\begin{thm}\label{thmOneg} Let $H={\rm SO}^-(n,2), n=2m$.
Then  $a_1,\cdots,b_4,\l_1,\l_2$ are
\begin{align*}
& a_1=0,&&
    a_2 = 2^{n-2},\\
&    a_3 = 2^{n-1} - 2^{m-1} - 1,&&
    a_4 = 2^{n-2} - 2^{m-1}  - 2,\\
 &   b_1 = 2^{n-3} + 2^{m-2} ,&&
    b_2=0,\\
    &b_3 = 2^{n-1} + 2^{m-1} ,&&
    b_4 = 2^{n-2} -1,\\
  &  \l_1 = 2^{n-2}- 1,&&
    \l_2 =2^{n-2} + 2^{m-1} .\qed 
\end{align*}
\end{thm}

\section { Relations for the $r=3$ cases and their groups}

\subsection {The $S_{n+1}$ case for $r=3$}

Here we have $b_3=\binom {n+1} 2,   b_4=\binom {n-1}2.$ Using these values for $b_3,b_4$ we solve  equations (\ref {ir13})-(\ref {ir63}) 
to obtain $\l_2=N/D$, where
$$
N:=\frac 1 4 2^n n^2 + \frac 1 4 2^n n - \frac 1 4 n^4 + \frac 1 2 n^3 - \frac 5 4 n^2 - 2n; \quad D=2^n - \frac 1 2 n^2 - \frac 1 2 n - 3.
$$
Now 
$$
N-\left (\frac 1 4 n^2+\frac 1 4 n\right )D=    -\frac 1 8 n^4 + \frac 3 4 n^3 - \frac 5  8 n^2 -\frac  3 4 n.
$$
But $N,D, N/D \in \mathbb Z$ and $\frac 1 4 n^2+\frac 1 4 n \in \frac 1 2 \mathbb Z$ gives
$$
A:=\frac {  -\frac 1 8 n^4 + \frac 3 4 n^3 - \frac 5  8 n^2 -\frac  3 4 n } { 2^n - \frac 1 2 n^2 - \frac 1 2 n - 2} \in \frac 1 2 \mathbb Z,
$$
and one checks that $A$ is not a half-integer for $n\ge 6$, a contradiction.

\subsection {The $S_{n+2}$ case for $r=3$}
Here we have $b_3=\binom {n+2} 2,   b_4=\binom {n}2.$ Using these values for $b_3,b_4$ we solve  equations (\ref {ir13})-(\ref {ir63}) 
to obtain $\l_2=N/D$, where
$$
N= 2^{n-2} n^2 + \frac 3 4  2^n n +2^{n-1} -\frac 1 4 n^4 - \frac 1 2 n^3 -\frac 5 4 n^2-4n-3;\quad D=2^n-\frac 1 2 n^2-\frac 3 2 n-3.
$$ 
Then 
$$
N-\left (\frac 1 4 n^2+\frac 3 4 n+\frac 1 2\right )D=-\frac 1 8 n^4+\frac 1 4 n^3+\frac 7 8 n^2- n-\frac 3 2.
$$
Then $N,D, N/D \in \mathbb Z$ and $\frac  1 4 n^2+\frac 3 4 n+\frac1 2      \in \frac 1 2 \mathbb Z$ gives
$$
\frac {-\frac 1 8 n^4+\frac 1 4 n^3+\frac 7 8 n^2- n-\frac 3 2 }{  2^n-\frac 1 2 n^2-\frac 3 2 n-4} \in \frac 1 2 \mathbb Z.
$$
One checks that the left-hand-side is not a half-integer for $n\ge 6$, a contradiction.

\subsection {The ${\rm SO}^+(n,2)$ case for $r=3$}
Let $n=2m$.
Here we have $b_3=2^{m-1}(2^m-1), b_4=2^{n-2}-1$. Using these values for $b_3,b_4$ we solve  equations (\ref {ir13})-(\ref {ir63}) 
to obtain $\l_2=N/D$, where
$$
N=\frac 1 4 2^{n} -\frac  1 4 2^{n+m} - 2^{n+1} + 2^{m+1};
\quad D=
2^n + 2^m - 4.
$$
Then $N-(2^{n-2}-2^{m-1}-1/2)D=2^{m-1}-2$ and 
 one checks that $\frac {  2^{m-1}-2 } {   2^n + 2^m - 4 }$  is not a half-integer if $n\ge 4$.

 \subsection {The ${\rm SO}^-(n,2)$ case for $r=3$}  
Let $n=2m$.
Here we have $b_3=2^{m-1}(2^m+1), b_4=2^{n-2}-1$. 
Using these values for $b_3,b_4$ we solve  equations (\ref {ir13})-(\ref {ir63}) 
to obtain $\l_2=N/D$, where
$$
N=\frac 1 4 2^{2n} + \frac 1 4 2^{n+m} - 2^{n+1} - 2^{m+1};\quad D=2^n - 2^m - 4.
$$
Then $N-(2^{n-2}+2^{m-1}-1/2)D=-2^{m-1}-2$ and 
 one checks that $\frac {  -2^{m-1}-2 } {   2^n - 2^m - 4 }$  is not a half-integer if $n\ge 4$.
 
 %%%%%%%%%%% use Spn2paperidealcalcs file
 
 \medskip
 
 Having checked these four cases we conclude:
 
 \begin{prop}\label{thm33333} The case $r=3$ cannot happen if $n\ge 6$.\qed
 \end{prop} 
 
 So we now assume that $r=2$ and $H=\langle C_2\rangle={\rm SO^\varepsilon}(n,2), \varepsilon =\pm$.

\section { The $r=2, H={\rm SO}^\pm(n,2) \le {\rm Sp}(n,2)$ cases}

The following result will allow us to reduce to the $H={\rm SO}^+(6,2)$ case.
Let $P_k$ denote the path graph with $k$ vertices.
\begin{prop} \label{thmspOn}
For  $n\ge 10$ even and $\varepsilon \in \{+,-\}$, the subgroup ${\rm SO}^\varepsilon(n,2)\le {\rm Sp}(n,2)$ contains a subgroup isomorphic to ${\rm SO}^+(6,2)\cong S_8$. 
\end{prop}
\noindent
{\it Proof} Let $n\ge 10$ be even and let $\e \in \{+,-\}$. 
We choose a basis $B$ using Figure 1 so that $\langle B \rangle={\rm SO}^\e(n,2)$ and $B$ contains vertices labeled $1,\cdots, n$, where
vertex $n$ is labeled $O^\varepsilon$ in Figure 1.
 Then referring to Figure 1 we see that there is a subset $C=\{b_1,b_2,b_3,b_4,b_5,b_6\}  \subset B$ such that $\Gamma^1(\{b_1\},\{b_2\},\{b_3\},\{b_4\},\{b_5\},\{b_6\})$ is the path  graph $P_6$.
 Then  $\Gamma^1(\{b_1\},\{b_2\},\{b_3\},\{b_4\},\{b_5\},\{b_6\},\{b_1+b_3+b_5\}\cong P_7$ and     $\langle C \cup \{b_1+b_3+b_5\}\rangle \cong S_8={\rm SO}^+(6,2)$.
\qed\medskip

Now consider the case where $r=2, H={\rm SO}^\pm(n,2)$ and let $b_1,b_2,\cdots,b_6$ be a part of a basis $B$ as in   Proposition \ref {thmspOn}  and Figure 1 where $\langle B\rangle =H$ and where 
$  t_{b_1},t_{b_2},\cdots ,t_{b_6}, t_{b_1+b_3+b_5} \in H$ generate $S_8\cong {\rm SO}^+(6,2)\le {\rm SO}^\pm (n,2)$.
Let $W:={\rm Span}(b_1,\cdots,b_6)$.
Then $ t_{b_1},t_{b_2},\cdots ,t_{b_6}, t_{b_1+b_3+b_5}  \in C_2 \cap \mathcal T(W)$ and the elements of $C_1 \cap \mathcal T(W)$ include those elements 
$e\in W$ where $Q_B(e)=0$.
We look at $(C_1  \cap \mathcal T(W))^3$, noting that it has multiplicities $1, 3, 4, 6, 20, 87$ and we let $X_i$ be those elements of the set $(C_1  \cap \mathcal T(W))^3$ that have multiplicity $i$. We show $X_1X_3 \ne X_3X_1$.

Define $\a=\a_1\a_2$ where 
$$
\a_1=t_{b_1+b_2+b_4}t_{b_2+b_5}t_{b_1+b_5}\in X_1,\, \a_2=t_{b_1+b_2+b_3+b_4+b_5}t_{b_2+b_3+b_5+b_6}t_{b_1+b_2+b_3+b_5+b_6}\in X_3,
$$
and note that the span of 
 $b_1+b_2+b_4,b_2+b_5,b_1+b_5,  
b_1+b_2+b_3+b_4+b_5,b_2+b_3+b_5+b_6,b_1+b_2+b_3+b_5+b_6$
 is $W$. This shows (i) that the minimal number of transvections needed to give $\a$ is six, and (ii) that
if $\a=t_ut_vt_wt_xt_yt_z$, then $u,v,w,x,y,z \in W$.
Now the multiplicities of $\a$ in $(X_1\cap \mathcal T(W)^3)(X_3\cap \mathcal T(W)^3)$ and $(X_3\cap \mathcal T(W)^3)(X_1\cap \mathcal T(W)^3)$ are $1$ and $3$ (respectively) and so  
$(X_1\cap \mathcal T(W)^3)(X_3\cap \mathcal T(W)^3)\ne (X_3\cap \mathcal T(W)^3)(X_1\cap \mathcal T(W)^3)$. 
It follows from (i) and (ii) that $X_1X_3 \ne X_3X_1.$

This concludes consideration of all cases and so the proof of Theorem \ref {mainthm1}.\qed

 \medskip
 
\noindent {\bf Proof that ${\rm Dye}_n, n>6,$  is not a strong Gelfand subgroup:}
We note from \cite {D1},\cite {D2} that if  $n=2m$ with $m$ prime, then  ${\rm Dye}_n\cong {\rm SL}(2,2^m)\rtimes \mathcal C_m$. 
We use a result from \cite [Lemma 3.3] {GH}:  For a group $H$ with irreducible character degrees $d_1,\cdots,d_t$ the {\it total character degree of $H$} is
$\tau_H:=\sum_{i=1}^t d_i$;  the result is: if $H \le G$ and $\chi \in \hat G,$ where $\tau_H< \chi(1)$, then $(G,H)$ is not a strong Gelfand pair. We apply this to ${\rm Dye}_n\le G_n$, where  $\tau_{{\rm Dye}_n}<|{\rm Dye}_n|=(2^n-1)2^mm$ and where $\chi\in \hat {G_n}$ is the Steinberg   character so that  $\chi(1)=2^{m^2}$, the size of a Sylow $2$-subgroup of $G_n$. Since 
$(2^{2m}-1)2^mm<2^{m^2}$ for $m>3$, the result follows.\qed

 \end{document}